\documentclass{amsart}

\usepackage{amsmath,amsthm,amsfonts,amssymb,bm,graphicx, mathrsfs}

\usepackage
{hyperref}
\hypersetup{colorlinks=true,citecolor=blue,linkcolor=blue,urlcolor=blue,
pdfstartview=FitH }

\theoremstyle{plain}
\newtheorem{theorem}{Theorem}

\newtheorem{lemma}{Lemma}

\numberwithin{equation}{section}

\theoremstyle{definition}

\renewcommand{\geq}{\geqslant}
\renewcommand{\leq}{\leqslant}

\newcommand{\changed}[1]{{\color{red} #1}}

\usepackage{calc}
\newsavebox\CBox
\newcommand\hcancel[2][0.5pt]{%
  \changed{\ifmmode\sbox\CBox{$#2$}\else\sbox\CBox{#2}\fi%
  \makebox[0pt][l]{\usebox\CBox}%
  \rule[0.5\ht\CBox-#1/2]{\wd\CBox}{#1}}}
\newcommand{\Oneb}{2/595}

\newcommand{\Halfbfrac}{\frac{1}{595}}

\newcommand{\Eighthbfrac}{\frac{1}{2380}}
\newcommand{\deltaFinal}{1/35000}
\newcommand{\deltaQuarter}{1/140000}
\newcommand{\deltaQuarterfrac}{\frac{1}{140000}}

\newcommand{\Deltavar}{\Delta}
\newcommand{\Xvar}{X}

\makeatletter
\DeclareRobustCommand\widecheck[1]{{\mathpalette\@widecheck{#1}}}
\def\@widecheck#1#2{%
    \setbox\z@\hbox{\m@th$#1#2$}%
    \setbox\tw@\hbox{\m@th$#1%
       \widehat{%
          \vrule\@width\z@\@height\ht\z@
          \vrule\@height\z@\@width\wd\z@}$}%
    \dp\tw@-\ht\z@
    \@tempdima\ht\z@ \advance\@tempdima2\ht\tw@ \divide\@tempdima\thr@@
    \setbox\tw@\hbox{%
       \raise\@tempdima\hbox{\scalebox{1}[-1]{\lower\@tempdima\box
\tw@}}}%
    {\ooalign{\box\tw@ \cr \box\z@}}}
\makeatother

\begin{document}

\author{Valentin Blomer}

\author{Jack Buttcane}
\address{Mathematisches Institut, Bunsenstr. 3-5, 37073 G\"ottingen, Germany} \email{vblomer@math.uni-goettingen.de}
\address{224 Mathematics Bldg, Buffalo, NY 14260, USA}
\email{buttcane@buffalo.edu}
 
 \title{Subconvexity for $L$-functions of non-spherical cusp forms on ${\rm GL}(3)$}

\thanks{First author partially supported by  DFG grant BL 915/2-2. Second author supported by NSF grant DMS-1601919}


\begin{abstract} Let $f$ be a cusp form for  ${\rm SL}(3, \Bbb{Z})$ associated with a generalized principal series representation of minimal weight $d$,  spectral parameter $r$ and associated $L$-function $L(s, f)$.  
For   $r \asymp d \asymp T$   the subconvexity bound $L(1/2, f) \ll  T^{\frac{3}{4} - \deltaQuarterfrac}$ is proved. 
\end{abstract}


\setcounter{tocdepth}{2}  \maketitle 

\maketitle

\section{Introduction} The two most commonly studied types of cusp forms for ${\rm GL}_3(\Bbb{A}_{\Bbb{Q}})$ are Maa{\ss} forms that are spherical at infinity, i.e.\ right ${\rm SO}(3)$-invariant, or symmetric square lifts of holomorphic forms of some weight $k$ for a congruence subgroup of ${\rm SL}_2(\Bbb{Z})$. 
There is a    family of non-spherical cusp forms that  is well-understood from the representation theoretic point of view and furnishes the so-called generalized principal series.  At the infinite place it is parametrized by a pair $(d, r)$, where  $d \geq 2$ is an integer, the (minimal) weight coming from an induced discrete series representation of the $2$-by-$2$ block of the $(2, 1)$-Levi subgroup of ${\rm GL}(3)$, and $r\in i\Bbb{R}$ is the spectral parameter. In technical terms, the underlying irreducible unitary representation of ${\rm GL}_3(\Bbb{R})$ is given as a Langlands quotient $\pi_{\infty} = J({\rm GL}_3(\Bbb{R}), P_{2, 1}, \sigma_1[r], \sigma_2)$ where $\sigma_1$ is a discrete series representation of weight $d$ and $\sigma_2(x) = |x|^{-2r}$ or $\text{sgn}(x) |x|^{-2r}$ depending on whether $d$ is even or odd. 

Cusp forms associated with the generalized principal series  transform with respect to ${\rm SO}(3)$ in terms of the $(2d+1)$-dimensional representation; thus when interpreted as functions on the generalized upper half plane $\Bbb{H}_3 = {\rm SL}_3(\Bbb{R})/{\rm SO}(3)$, they come as $(2d+1)$-dimensional vector valued automorphic forms. Symmetric square lifts for holomorphic forms of weight $k$ correspond to forms with $(d, r) = (2k-1, 0)$. 

It was proved only recently by the second author \cite{Bu2} that cusp forms associated with the generalized principal series (that are not symmetric squares of holomorphic forms) exist in abundance: their spectral density is 
\begin{equation}\label{spec}
   \textbf{spec}^d(r) = \frac{1}{8 \pi^3 } (d-1) \left(\frac{1}{4}(d-1)^2 -9r^2\right) \asymp d(d^2 + |r|^2)
   \end{equation} and the weighted Weyl law \cite[Theorem 1]{Bu2} together with \eqref{l1ad} below shows that for $D, R$ sufficiently large there are   $\gg_{\varepsilon} (D^2R(D^2 + R^2))^{1-\varepsilon}$ such forms with  $d \asymp D$, $|r| \asymp R$.\\ 
 
 While there exists extensive literature on spherical forms and symmetric square lifts, we would like to demonstrate in action   the analytic theory of typical  Maa{\ss} forms belonging   to the generalized principal series.  We investigate the $L$-functions $L(s, f)$ of such cusp forms $f$, whose Dirichlet series feature the usual Hecke eigenvalues $A_f(1, n)$ and whose gamma factors  
  \begin{equation}\label{linf}
 L_{\infty}(s, f) :=  \Gamma_{\Bbb{C}}(\textstyle\frac{1}{2}(d-1)\displaystyle + r + s) \Gamma_{\Bbb{R}}(\mathfrak{a} + s - 2r), \quad \{0, 1\} \ni \mathfrak{a} \equiv d \, (\text{{\rm mod }} 2)
 \end{equation}
with $\Gamma_{\Bbb{R}}(s) = \pi^{-s/2} \Gamma(s/2)$, $\Gamma_{\Bbb{C}}(s) = 2(2\pi)^{-s} \Gamma(s)$ feature the spectral data $(d, r)$, see \cite[p.\ 101]{HIM} or \cite[A.3]{RS}. 
 In particular, the analytic conductor $\mathcal{C}(f)$ of $f$ at the central point $s = 1/2$ is of size 
 \begin{equation}\label{conductor}
   \mathcal{C}(f) \asymp (1+|r|)(d^2 + |r|^2),
 \end{equation}  
    and therefore the  convexity bound is 
 \begin{equation}\label{conv}
    L(1/2, f) \ll \left( (1+|r|) (d^2  + |r|)^2\right)^{1/4 + \varepsilon}.
 \end{equation}   
     If $T > 1$ is a large parameter, a generic form $f$ of conductor $\mathcal{C}(f) \asymp T^3$ will have 
 \begin{equation}\label{rdT}
    |r| \asymp d \asymp T,
    \end{equation}  and we will focus on such forms, which cover 99.9\% of all forms with conductor of size $T^3$.      Our main result is the following subconvexity bound, improving on the bound \eqref{conv}.  
 \begin{theorem}\label{mainthm} Let $f$ be a cusp form of ${\rm GL}(3)$ over $\Bbb{Q}$ that is unramified at all finite places and associated to the generalized principal series of weight $d$ and spectral parameter $r$. Suppose that $T > 1$ is sufficiently large and \eqref{rdT} holds. Then
 \begin{equation}\label{mainbound}
 L(1/2, f) \ll T^{3/4 - \deltaQuarter}.
 \end{equation}
 \end{theorem}
 This is the non-spherical analogue to the subconvexity result in \cite{BB}, and we assume some familiarity with the method of  \cite{BB}. As in \cite{BB} we have made no effort to optimize the exponent. We have incorporated a number of simplifications and improvements, some of which are specific to the non-spherical case under consideration.  The situation is roughly comparable to, but more complicated than,  the analytic features of the Petersson formula versus the Kuznetsov formula: the Bessel $J$-function with \emph{real} index decays very quickly for arguments less than the index, but on the other hand  the Bessel $J$-function with \emph{imaginary} index has no transitional range (for real arguments). Correspondingly, Lemmas \ref{lem3} and \ref{lem4} are simpler than the corresponding  Lemmas 8 and 9 in \cite{BB}, but  parts of Lemma \ref{lem6} require more work than Lemma 11 in \cite{BB}, and this is responsible for the numerically slightly weaker exponent. The  bound \eqref{mainbound} holds more generally at any point $1/2 + it$, $t \in \Bbb{R}$ fixed, on the critical line.
  
 The basic idea is to estimate an amplified fourth moment using a version of the Kuznetsov formula for this particular spectral family. We apply Poisson summation on all four variables, which transforms the long Weyl element Kloosterman sums essentially into congruences. To solve the resulting counting problem, we have to understand the 4-fold Fourier transform of the corresponding Whittaker-type kernels in the Kuznetsov formula, which requires a very subtle analysis. \\
  
  We hope that this first analytic result for $L$-functions of generalized principal series will stimulate further research in this direction. Several variations on this theme are possible. For instance, one can relax \eqref{rdT} and require only $|r| \asymp T$, $d \asymp T^{\alpha}$ with $1/2 \leq \alpha \leq 1$. This does not change the size of the analytic conductor \eqref{conductor}, but shortens the family. One can recover the size $T^3$ of the family (and hence the structure of the present analysis) by averaging over forms of weight $[d, d + T^{1-\alpha}]$. These and other generalizations are left to future investigations.\\
  
A word on \emph{notation}: as usual, we use $\varepsilon$-convention, meaning that the letter $\varepsilon$ denotes a sufficiently small number whose value may change at each occurrence. For non-zero, but not necessarily positive real quantities $A, B$ we write $A \asymp B$ to mean that there are positive constants $c_1, c_2$ such that $c_1 A \leq B \leq c_2A$. The word negligible is synonymous to $\ll_B T^{-B}$ for any $B \geq 0$.

 \section{The Kuznetsov formula}
 
The statement of the Kuznetsov formula for non-spherical forms requires  quite a bit of notational preparation. We follow the notation of \cite{Bu2} and \cite{BB}. 

\subsection{The spectral side}  We start by observing that cusp forms for the generalized principal series satisfy automatically the generalized Ramanujan conjecture at infinity.  This follows from the fact that the Langlands parameters $\mu$ are given by $$\mu  = (\textstyle\frac{1}{2}(d-1)+ r, -\frac{1}{2}(d-1) + r, -2r)$$ together with the unitarity condition that the entries $\mu$ are a permutation of the entries of $-\bar{\mu}$, so that necessarily $r\in i\Bbb{R}$. 

For $d \in \{2, 3, 4\ldots\}$, $r \in i\Bbb{R}$, $0 \leq m \leq d$, $\epsilon \in \{\pm 1\}$, $y  = (y_1, y_2) \in \Bbb{R}_{>0}^2$ we define $W_{\epsilon m}(y; d, r)$ to be 
\begin{displaymath}
\begin{split}
  & \frac{1}{2^{1+d} \pi} 
\left(\begin{matrix} 2d\\ d+m \end{matrix}\right)^{1/2} \sum_{\ell = 0}^m \epsilon^{\ell} \left(\begin{matrix}m \\ \ell \end{matrix}\right) \int_{(c)} (2\pi y_1)^{1-s_1} (2\pi y_2)^{1-s_2} \\
&\times \Gamma \left(\frac{d-1}{2} + s_1 - r\right) \Gamma \left(\frac{d-1}{2} + s_2 + r\right) B\left(\frac{d-m+s_1+2r}{2}, \frac{\ell + s_2 - 2r}{2}\right)  \frac{ds}{(2\pi i)^2}
\end{split}
\end{displaymath}
where $c = (c_1, c_2)$ is a pair of positive real numbers and $B$ denotes the beta function.  The  completed Whittaker function  is the row vector
$$W(y; d, r) = \big(W_{-d}(y; d, r), \ldots, W_d(y; d, r)\big).$$
We use the usual coordinates on the generalized upper half plane $\Bbb{H}_3 $ 
and write 
$x= \left(\begin{smallmatrix} 1& x_2 & x_3\\   & 1 & x_1\\   &  & 1\end{smallmatrix}\right)$, $ y = \text{diag}(y_1 y_2, y_1, 1)$
 with measures $dx = dx_1\, dx_2\, dx_3$, $d^{\ast}y = dy_1\, dy_2/(y_1y_2)^3$. 


If $f$ is a cusp form of weight $d$ and spectral parameter $r$, viewed as a $(2d+1)$-dimensional vector valued form $(f_{-d}, \ldots, f_d)$ on ${\rm SL}_3(\Bbb{Z}) \backslash {\rm SL}_3(\Bbb{R})$, we normalize its Fourier coefficients $\rho_f(n_1, n_2)$, $n_1, n_2 \in \Bbb{Z} \setminus \{0\}$, by setting
$$\int_{U(\Bbb{Z}) \backslash U(\Bbb{R})} f\left(\left(\begin{smallmatrix} 1& x_2 & x_3\\   & 1 & x_1\\   &  & 1\end{smallmatrix}\right)\left(\begin{smallmatrix} y_1y_2 &   & \\   & y_1 &  \\   &  & 1\end{smallmatrix}\right)\right) e(-x_1n_1 - x_2n_2) dx =  \frac{\rho_f(n_1, n_2)}{|n_1n_2|} W\left( (n_1y_1, n_2y_2) ; d, r\right),$$
where as usual $U$ denotes the algebraic group of unipotent upper triangular matrices. 
The inner product is given by $$\langle f, g \rangle =  \int_{{\rm SL}_3(\Bbb{Z}) \backslash \Bbb{H}_3} f(x y) \overline{f(x y)}^{\top}  dx\, d^{\ast}y.$$ 

We write
\begin{displaymath}
\begin{split}
&  \int_{[0, \infty)^2} W(y; d, r) W(y; d, -r)^{\top}  y_1^2y_2  \, d^{\ast}y =: \frac{1}{\textbf{cos}^d(r)}.
\end{split}
\end{displaymath}
This and more general expressions are computed explicitly in \cite[Theorem 2]{Bu2} as a generalization of Stade's formula. 
 If $f$ is a Hecke eigenform, then $\rho_f(n_1, n_2)$ is proportional to the Hecke eigenvalues $A_f(n_1, n_2)$. We can compute the $L^2$-norm of a Hecke-normalized cusp form by  the usual Rankin-Selberg unfolding method. The maximal Eisenstein series $E(z, s; \textbf{1})$ of  weight 0 twisted by the constant function has a simple pole with residue $2\pi/(3 \zeta(3))$ at $s=1$, so that (cf.\ e.g.\ \cite[(20)]{BBR})
\begin{displaymath}
  \| f \|^2 
  = \frac{3\zeta(3)}{2\pi}  \underset{s=1}{\text{res}} \langle f E(., s, \textbf{1}), f \rangle = \frac{3\zeta(3)}{2\pi}   \underset{s=1}{\text{res}} \frac{L(s, f \times \bar{f})}{\zeta(3s)} \frac{1}{\textbf{cos}^d(r)} = \frac{3 \,L(1, \text{ad}, f)}{2\pi\,\textbf{cos}^d(r)} . 
\end{displaymath} 
 We conclude that if $f$ is $L^2$-normalized, then
\begin{equation}\label{norm}
  \frac{\rho_f(n_1, n_2) \overline{\rho_f(m_1, m_2)}}{\textbf{cos}^d(r)} =  \frac{2\pi}{3} \frac{A_f(n_1, n_2) \overline{A_f(m_1, m_2)}}{L(1, \text{ad}, f)}
\end{equation}
for $n_1, n_2, m_1, m_2 \in \Bbb{N}$. 
 
Eisenstein series exist only if $d$ is even, and it is only the maximal  parabolic Eisenstein series that come up in the spectral decomposition.  They are parametrized by an orthonormal basis  of holomorphic cusp forms $\varphi$ of weight $d$ for ${\rm SL}_2(\Bbb{Z})$ and a spectral parameter $r \in i\Bbb{R}$. We normalize the Fourier expansion in the same way as we did for cusp forms and denote the Fourier coefficients by $\rho_{\varphi, r}(n_1, n_2)$. They are proportional to Hecke eigenvalues, and the size of $\rho_{\varphi, r}(1, 1)$ is computed explicitly in \cite[Section 10]{Bu2}, but we do not need this for our purposes. 

\subsection{The arithmetic side} We define the usual ${\rm GL}(3)$ Kloosterman sums for the $w_4 = \left(\begin{smallmatrix}   &1& \\ & & 1\\ 1&  &\end{smallmatrix}\right)$,  $w_5 = \left(\begin{smallmatrix}   && 1\\ 1 & & \\ & 1 &\end{smallmatrix}\right)$ and  $w_6 = -\left(\begin{smallmatrix}   &&1 \\ & 1& \\ 1 &  &\end{smallmatrix}\right)$  Weyl elements. For $n_1, n_2, m_1, m_2, D_1, D_2 \in \Bbb{N}$ we write 
\begin{displaymath}
	\tilde{S}(n_1,n_2,m_1;D_1,D_2) := \sum_{\substack{C_1 (\text{mod }D_1), C_2 (\text{mod }D_2)\\(C_1,D_1)=(C_2,D_2/D_1)=1}} e\left(n_2\frac{\bar{C_1}C_2}{D_1}+m_1\frac{\bar{C_2}}{D_2/D_1}+n_1\frac{C_1}{D_1}\right)
\end{displaymath}
if $D_1\mid D_2$,   and
\begin{displaymath}
\begin{split}
&S(n_1, m_2, m_1, n_2; D_1, D_2)\\
& := \sum_{\substack{B_1, C_1 \, ({\rm mod }\, D_1)\\B_2, C_2 \, ({\rm mod }\,  D_2)\\ D_1C_2 + B_1B_2 + D_2C_1 \equiv 0 \, ({\rm mod }\, D_1D_2)\\ (B_j, C_j, D_j) = 1}}\hspace{-0.9cm} e\left(\frac{n_1B_1 + m_1(Y_1 D_2 - Z_1 B_2)}{D_1} + \frac{m_2B_2 + n_2(Y_2 D_1 - Z_2B_1)}{D_2}\right),
\end{split}
\end{displaymath}
where the numbers $Y_j$ and $Z_j$ are defined by $Y_jB_j + Z_jC_j \equiv 1 \, (\text{mod }D_j)$ for $j = 1, 2$.

We define
\begin{equation}\label{qds}
Q(d, s) := \frac{\Gamma(\frac{d-1}{2} + s)}{\Gamma(\frac{d+1}{2} - s)}.
\end{equation}
For $y \in \Bbb{R} \setminus \{0\}$, $\epsilon = \text{sgn}(y)$, $d \in \Bbb{N}$, $d \geq 2$, $r \in i\Bbb{R}$,  we define  
\begin{equation}\label{kw4}
\begin{split}
K_{w_4}(y; d, r) & = \frac{(\epsilon i)^d}{4\pi^2} \int_{-i\infty}^{i \infty} |8\pi^3 y|^{1-r-s} Q(d, s) \Gamma(s+3r) \exp\left(\frac{\epsilon i \pi}{2}(s + 3r)\right) \frac{ds}{2\pi i}\\
\end{split}
\end{equation}
with the usual Barnes convention that the path of integration should pass to the right of all of the poles of the gamma functions in the form $\Gamma(s_j + a)$. 
Moreover, we choose the contour such that all integrals are absolutely convergent, which can always be arranged by shifting the unbounded part appropriately.

For $y  = (y_1, y_2) \in (\Bbb{R} \setminus \{0\})^2$, $\epsilon  = (\text{sgn}(y_1), \text{sgn}(y_2))$, we define
\begin{equation}\label{kw6}
\begin{split}
&K_{w_6}(y; d, r) \\
&= \frac{1}{4\pi^2} \left| \frac{y_2}{y_1}\right|^r  \int_{-i\infty}^{i\infty} \int_{-i\infty}^{i\infty} |4\pi^2 y_1|^{1-s_1} |4\pi^2 y_2|^{1-s_2} B_{w_6}^{\epsilon}((s_1, s_2), r) Q(d, s_1) Q(d, s_2) \frac{ds_2\, ds_1}{(2\pi i)^2}
\end{split}
\end{equation}
where
\begin{equation}\label{B}
B_{w_6}^{\epsilon}((s_1, s_2), r) = \begin{cases} (-1)^d B(s_1 + 3r, s_2 - 3r), & \epsilon = (-, -), \\ B(s_2-3r, 1 - s_1 - s_2), & \epsilon = (-, +) , \\ B(s_1+3r, 1 - s_1 - s_2), & \epsilon= (+, -), \\ 0, & \epsilon = (+,+). \end{cases}
\end{equation}

For a Schwartz class function $F$ that is holomorphic in $|\Im r| < 1/4 + \delta$ for some $\delta > 0$, we define
\begin{equation}\label{phi}
\begin{split}
& \Phi_{w_4}(y; d) := \frac{1}{|y|} \int_{(0)} F(r) K_{w_4}(y; d, r) \textbf{spec}^d(r)\, \frac{dr}{2\pi i},  \quad y \in \Bbb{R} \setminus \{0\},\\
& \Phi_{w_5}(y; d) := \frac{1}{|y|} \int_{(0)} F(r) K_{w_4}(-y; d, -r) \textbf{spec}^d(r)\, \frac{dr}{2\pi i}, \quad y \in \Bbb{R} \setminus \{0\},\\
& \Phi_{w_6}(y; d) := \frac{1}{|y_1y_2|} \int_{(0)} F(r) K_{w_6}(y; d, r) \textbf{spec}^d(r)\, \frac{dr}{2\pi i}, \quad y = (y_1, y_2) \in (\Bbb{R} \setminus \{0\})^2
\end{split}
\end{equation}
with $\textbf{spec}^d(r)$ as in \eqref{spec}.

\subsection{The Kuznetsov formula}
We are now ready to state the formula. Let $F$ be a Schwartz class function that is holomorphic in $|\Im r| < 1/4 + \delta$ for some $\delta > 0$, and keep the notation from the previous two subsections. Let $S_3^d$ be an orthonormal basis of cuspidal Hecke eigenforms associated with the generalized principal series of weight $d$, and for $f \in S_3^d$ let $r_{f} \in i\Bbb{R}$ be the spectral parameter. Let $S_2^d$ denote an orthonormal  basis of holomorphic cusp forms of weight $d$ (empty when $d$ is odd). Let $n_1, n_2, m_1, m_2 \in \Bbb{N}$. Let $d \geq 2$ be an integer. Then we have the following absolutely convergent summation formula \cite[Theorem 4]{Bu2}: 
\begin{equation}\label{kuz}
\begin{split}
& \sum_{f \in S_3^d} F(r_f) \frac{\overline{\rho_f(m_1, m_2)}\rho_f(n_1, n_2)}{\textbf{cos}^d(r_f)}  + 2\sum_{\varphi \in S_2^d}  \int_{(0)} F(r) \frac{\overline{\rho_{\varphi, r}(m_1, m_2)}\rho_{\varphi, r}(n_1, n_2)}{\textbf{cos}^d(r)} \frac{dr}{2\pi i}\\
&= \Delta^d + \Sigma^d_4 + \Sigma^d_5 + \Sigma^d_6
\end{split}
\end{equation}
where
\begin{displaymath}
\begin{split}
  \Delta^d& = \delta_{n_1, m_1} \delta_{n_2, m_2}  \int_{(0)} F(r) \textbf{spec}^d(r) \frac{dr}{2\pi i},\\
   \Sigma^d_{4}& = \sum_{\epsilon  = \pm 1} \sum_{\substack{D_2 \mid D_1\\  m_2 D_1= n_1 D_2^2}}\frac{ \tilde{S}(-\epsilon n_2, m_2, m_1; D_2, D_1)}{D_1D_2} \Phi_{w_4}\left(  \frac{\epsilon m_1m_2n_2}{D_1 D_2}; d \right),  \\ 
 \Sigma^d_{5} &= \sum_{\epsilon  = \pm 1} \sum_{\substack{     D_1 \mid D_2\\ m_1 D_2 = n_2 D_1^2}} \frac{ \tilde{S}(\epsilon n_1, m_1, m_2; D_1, D_2) }{D_1D_2}\Phi_{w_5}\left( \frac{\epsilon n_1m_1m_2}{D_1 D_2}; d\right),\\
   \Sigma^d_6 &= \sum_{\epsilon_1, \epsilon_2 = \pm 1} \sum_{D_1,  D_2  } \frac{S(\epsilon_2 n_2, \epsilon_1 n_1, m_1, m_2; D_1, D_2)}{D_1D_2} \Phi_{w_6}  \left( \left( \frac{\epsilon_2 m_1n_2D_2}{D_1^2},  \frac{\epsilon_1 m_2n_1D_1}{ D_2^2} \right); d\right). 
\end{split}
\end{displaymath} 

Given a large parameter $T > 1$ and a small $\varepsilon > 0$, as well as the spectral parameter $r_f \in i\Bbb{R}$ of our favourite Maa{\ss} form $f$ with $|r_f| \asymp T$, we will choose 
\begin{equation}\label{defF}
F(r) := \exp\left(\frac{(r-r_f)^2}{T^{\varepsilon}}\right) ,
\end{equation}
which satisfies the conditions of the Kuznetsov formula. Note that by \eqref{spec} we have 
\begin{equation}\label{intF}
\int_{(0)}  F(r)  \textbf{spec}^d(r) \frac{dr}{2\pi i} \ll T^{3+\varepsilon}
\end{equation}
for $d \asymp T$. 

 \section{Special functions}

We need the following well-known results for the Bessel $J$-function for $d \in \Bbb{N}$, $y > 0$. We start with the bound
\begin{equation}\label{jbound}
  J_d(y) \ll \min( (2y/d)^d, 1). 
\end{equation} 
The first bound follows from \cite[Lemma 4.1]{Ra} for $d \geq 15$  and from the power series expansion \cite[8.402]{GR}  for $1 \leq d \leq 14$ and  $x \leq 1$, the second bound follows from \cite[8.411.1]{GR}. The second bound can be regarded as a trivial bound (and could be improved, but this would have no influence for the present paper), while the first bound tells us that $J_d(y)$ essentially disappears if $y \leq d/3$, say, and $d$ is large. Again this can be improved, but this is irrelevant for the present purpose. 
By \cite[8.471.2]{GR}, the derivatives satisfy
\begin{equation}\label{jder}
  J'_d(y) = \frac{1}{2} (J_{d-1}(y)  - J_{d+1}(y)).
\end{equation}
We have the Mellin formula \cite[17.43.16]{GR}
\begin{equation}\label{mellin}
J_{d-1}(2\sqrt{x}) = \int_{(c)} Q(d, s)x^{-s}  \frac{ds}{2\pi i} 
\end{equation}
for $(1-d)/2 < c < 0$, $x > 0$ and $Q(d, s)$ as in \eqref{qds}.  We will us this in combination with the integral representation of the Euler beta function
\begin{equation}\label{beta}
B(x, y) = \frac{\Gamma(x)\Gamma(y)}{\Gamma(x+y)} = \int_0^1 t^{x-1}(1-t)^{y-1} dt , \quad \Re x, \Re y > 0,
\end{equation}
and the Mellin formula for the exponential function (\cite[17.43.3/4]{GR}) 
\begin{equation}\label{exp}
 \exp(\pm 2 i x) x^{3r}    =  \int_{-i\infty}^{i\infty} 2^{-3r-s} \Gamma(s + 3r) \exp\left( \pm \frac{1}{2}  i\pi (3r+s)\right) x^{-s} \frac{ds}{2\pi i}
 \end{equation}
for $r \in i\Bbb{R}$, $x > 0$. 

As in \cite[Section 5]{BB}, the kernel functions $K_{w_4}$ and $K_{w_6}$, defined in \eqref{kw4} and \eqref{kw6} as  Mellin-Barnes integrals, have alternative integral representations in terms of Bessel functions, and these turn out to be useful and of independent interest.  With $\epsilon = (\text{sgn}(y_1), \text{sgn}(y_2))$ we see  that   $K_{w_6}(y; d, r) $ equals
\begin{equation}\label{bessel}
 \begin{cases} \displaystyle (-1)^d 4\pi^3 \frac{|y_1y_2|}{|y_1/y_2|^r}  \int_0^1 x^{3r} (1-x)^{-3r} J_{d-1}\left(4\pi\sqrt{\frac{|y_1|}{x}}\right) J_{d-1}\left(4\pi\sqrt{\frac{|y_2|}{1-x}}\right)dx,&  \epsilon = (-, -),\\ \displaystyle  \frac{16\pi^4|y_1y_2|}{|y_1/y_2|^r} \int_0^1 x^{-3r} (1-x) J_{d-1}\left(4\pi\sqrt{ |y_1|(1-x)}\right) J_{d-1}\left(4\pi\sqrt{ |y_2|(x^{-1} - 1) }\right)dx, &  \epsilon = (-, +),\\
\displaystyle  \frac{16\pi^4|y_1y_2|}{|y_1/y_2|^r} \int_0^1 x^{3r} (1-x) J_{d-1}\left(4\pi\sqrt{ |y_1|(x^{-1} - 1)}\right) J_{d-1}\left(4\pi\sqrt{ |y_2|(1-x) }\right)dx, &  \epsilon = (-, +),\\
0, &\epsilon = (+, +). 
  \end{cases}
  \end{equation}
This can be checked easily by inserting the Mellin formula \eqref{mellin} for the Bessel functions and evaluating the remaining $x$-integral with \eqref{beta}. 

Similarly, $K_{w_4}(y; d, r)$ equals
\begin{equation}\label{bessel1}
   K_{w_4}(y; d, r) = 2(\epsilon i)^d  |y|^{1-r} \pi^{1-3r} \int_0^{\infty} J_{d-1}(2\sqrt{x}) \exp\left(2 \epsilon i \frac{4\pi^3|y|}{x}\right) \left(\frac{4\pi^3|y|}{x}\right)^{3r} \frac{dx}{x}.
\end{equation}
This follows by inserting   \eqref{exp} and evaluating the $x$-integral by \eqref{mellin} and Mellin inversion. 

\section{$L$-functions}
As before let $f$ be a cusp form associated with the generalized principal series of weight $d$ and spectral parameter $r = r_f \in i\Bbb{R}$. We assume that $f$ is a Hecke eigenform with Hecke eigenvalues $A_f(n, m)$ and denote by $$L(s, f) = \sum_{n} \frac{A_f(1, n)}{n^s}$$ the associated $L$-function. 
The completed version is $\Lambda(s, f) = L_{\infty}(s, f)L(s, f)$ with $L_{\infty}(s, f)$ as in \eqref{linf}. In particular, we can express the central value $L(1/2, f)$ by a standard approximate functional equation of length about $T^3$ if $d \asymp |r_f| \asymp T$. The Hecke eigenvalues (as well as those of Eisenstein series) satisfy the usual Hecke relations of the unramified ${\rm GL}(3)$ Hecke algebra  (\cite[Section 6]{Go}), in particular
$$A_f(1, n) \overline{A_f(1, m)} = \sum_{d \mid (n, m)} A_{f}\left(\frac{m}{d}, \frac{n}{d}\right).$$

 We follow the argument of \cite[Section 2]{BB}, so we can be brief. 
For a suitable  smooth compactly supported  function $W$ (all of whose derivatives are bounded independently of $T$) and $M \geq 1$ we define $$\mathcal{L}_M := \frac{1}{M} \Big| \sum_{n, m} A_{f}(m, n) W\left(\frac{n}{M} \right) \overline{W\left(\frac{m}{M}\right)}\Big|,$$
which (by Rankin-Selberg theory) satisfies the trivial bound 
\begin{equation}\label{etasmall}
\mathcal{L}_M \ll M(MT)^{\varepsilon}.
\end{equation}
As in \cite[(2.4)]{BB} we have
$$|L(1/2, f)|^2 \ll T^{\varepsilon} \max_{M \leq T^{3/2 + \varepsilon}} \mathcal{L}_M.$$
The trivial bound suffices if $M$ is small, so from now on we assume 
\begin{equation}\label{sizeM}
   T^{3/2 - \eta} \leq M \leq T^{3/2 + \varepsilon}
   \end{equation}
for some small $0 < \eta < 1$. Let $x(n) :=  A_{f}(1, n)/| A_{f}(1, n)|$ if $A_{f}(1, n) \not= 0$ and $x(n) = 0$ otherwise. Fix some sufficiently small $0 < \lambda < 1/20$, and for
\begin{equation}\label{lambda}
L = T^{\lambda} \leq T^{1/20}
\end{equation}
 define the amplifier 
$$\mathcal{A}(g) = \sum_{j=1}^3 \Bigl| \sum_{\substack{L \leq \ell \leq 2L\\ \ell \text{ prime}}} A_{g}(1, \ell^j) \overline{x(\ell^j)}\Bigr|^2$$
where $g$ can be a cuspidal Hecke eigenform in $S_3^d$ or an Eisenstein series associated to a pair $(\varphi, r)$ occurring on the spectral side of the Kuznetsov formula. The Hecke relation 
$$A_{f}(1, \ell) A_{f}(1, \ell^2) = A_{f}(1, \ell^3) + A_{f}(1, \ell)A_{f}(\ell, 1) - 1$$
implies $\mathcal{A}(f) \gg L^{2-\varepsilon}$ for our preferred form $f$.  The bound \cite[Theorem 2]{Li}
\begin{equation}\label{l1ad}
L(1, f, \text{ad}) \ll T^{\varepsilon},
\end{equation}
now implies 
$$\mathcal{L}_M^2  \ll  \frac{\mathcal{L}_M^2 \mathcal{A}(f) T^{\varepsilon}}{L^{2-\varepsilon} L(1, f, \text{ad}) }. $$
By \eqref{norm} and positivity\footnote{Which is why the normalization the exact shape of the Eisenstein series plays no role, we only need the abstract fact that they can be Hecke-diagonalized in order to use the amplifier at all spectral components.}  we can now take an appropriately weighted sum   over the spectrum  and   arrive   at the basic inequality
\begin{equation}\label{basic}
\begin{split}
\mathcal{L}_M^2  \ll&   \frac{T^{\varepsilon}}{M^2L^2}  \sum_{j=1}^3
 \sum_{\substack{\ell_1, \ell_2 \asymp L\\ \ell_1, \ell_2 \text{ prime}}}  \sum_{\substack{r_0r_1r_2 = \ell_1^j \\ s_0s_1s_2 = \ell_2^j }}   \Bigl|  \sum_{n_1, n_2, m_1, m_2} W\Bigl(\frac{r_2n_1}{M}\Bigr) W\Bigl(\frac{s_0m_1}{M}\Bigr)    \overline{W\Bigl(\frac{s_1n_2}{M}\Bigr)W \Bigl(\frac{r_0m_2}{M} \Bigr) }  \\
 & \quad\quad\times \Bigl\{\sum_{f \in S_3^d} F(r_f) \frac{\overline{\rho_f(m_1s_1, n_2s_2)}\rho_f(m_2r_2, n_1r_1)}{\textbf{cos}^d(r_f)}  \\
 & \quad\quad\quad\quad\quad\quad\quad + 2\sum_{\varphi \in S_2^d}  \int_{(0)} F(r) \frac{\overline{\rho_{\varphi, r}(m_1s_1, n_2s_2)}\rho_{\varphi, r}(m_2r_2, n_1r_1)}{\textbf{cos}^d(r)} \frac{dr}{2\pi i}\Bigr\}\Bigr| \\
 \end{split}
\end{equation}
with $F$ as in \eqref{defF}. 
This is in good shape to apply the Kuznetsov formula.

\section{Kloosterman sums}

In this section we quote two lemmas from \cite[Section 6]{BB}.

\begin{lemma}\label{lem1} Let $s_1, s_2, r_1, n_2 \in \Bbb{Z} \setminus \{0\}$. Let $D, \delta \in \Bbb{N}$ and $x, y \in \Bbb{Z}$. Then
\begin{equation}\label{Klo5}
\begin{split}
\Bigl|\frac{1}{D\delta}&\sum_{ n_1 \, (\text{{\rm mod }}D)}\sum_{ m_1 \, (\text{{\rm mod }} \delta)} \tilde{S}(n_1r_1, n_2s_2, m_1s_1; D, D\delta ) e\left( - \frac{xn_1}{D} - \frac{ym_1}{\delta}\right)\Bigr| \leq D (r_1,  D)(s_1,  \delta),  
\end{split}
\end{equation}
and the left hand side vanishes unless $(D, x) = (r_1, x)$,  $(\delta, y) = (s_1, y)$ and $\frac{D}{(D, \delta)} \mid n_2s_2$. 
\end{lemma}

For $r_1, r_2, s_1, s_2 \in \Bbb{Z} \setminus \{0\}$, $x_1, x_2, y_1, y_2 \in \Bbb{Z}$ and $D_1, D_2 \in \Bbb{N}$ we define the finite Fourier transform
\begin{equation}\label{hatS}
\begin{split}
&\widehat{S}_{r_1, s_1, r_2, s_2}(x_1, y_1, x_2, y_2; D_1, D_2) := \\
&\frac{1}{D_1^2D_2^2} \sum_{\substack{n_1, m_1\, (\text{{\rm mod }} D_1)\\ n_2, m_2 \, (\text{\rm{mod  }} D_2)}} S(n_1r_1, m_2r_2, m_1s_1, n_2s_2; D_1, D_2)e\left( - \frac{x_1m_1 + y_1 n_1}{D_1} - \frac{x_2m_2 + y_2 n_2}{D_2}\right).
\end{split}
\end{equation}
As usual we denote Euler's function by $\phi$.  For a prime $\ell$ we write $r \mid \ell^{\infty}$ if $r$ is a power of $\ell$, and we denote by $(\ell^{\infty}, r)$ the highest power of $\ell$ dividing $r$.

\begin{lemma}\label{lem2} {\rm (a)} We have the general bound
 \begin{equation*}\label{Klo6}
 \begin{split}
|\widehat{S}_{r_1, s_1, r_2, s_2}(x_1, y_1, x_2, y_2; D_1, D_2)| \leq  (r_1, D_1)  (r_2, D_2)  (D_1, D_2),  
  \end{split}
  \end{equation*}
  and the left hand side vanishes unless $x_1y_1 \equiv r_1s_1 D_2 \, (\text{{\rm mod }}D_1)$ and $x_2y_2 \equiv r_2s_2 D_1 \, (\text{{\rm mod  }}D_2)$. \\
  {\rm (b)}  We have $\widehat{S}_{r_1, s_1, r_2, s_2}(0, 0, x_2, y_2; D_1, D_2)  = 0$  unless
  $$\frac{D_1}{(D_1, r_1s_1)} \mid (x_2, y_2).$$
     Similarly,  
  $\widehat{S}_{r_1, s_1, r_2, s_2}(x_1, y_1, 0, 0; D_1, D_2)  = 0$ unless $\frac{D_2}{(D_2, r_2s_2)} \mid (x_1, y_1).$\\
  {\rm (c)} If $(r_1r_2, s_1s_2) = 1$, then $\widehat{S}_{r_1, s_1, r_2, s_2}(0, 0, 0, 0; D_1, D_2)  = 0$ unless $D_1 = D_2$, in which case it equals $\phi(D_1) = \phi(D_2)$. \\
  {\rm (d)} Let $\ell$ be a prime and assume that $r_1r_2s_1s_2 \mid \ell^{\infty}$. Then
 \begin{equation}\label{partd}
 |\widehat{S}_{r_1, s_1, r_2, s_2}(0, 0, 0, 0; D_1, D_2)| \leq (D_1, D_2) (\ell^{\infty}, [D_1, D_2])r_2.
 \end{equation}
   \end{lemma}
    

\section{Analysis of the Bessel kernels I} 

The heart of the analysis is contained in the  investigation of the  properties of the weight functions $\Phi_w$, defined in \eqref{phi} with $F$ as in \eqref{defF}.  The following three lemmas correspond to \cite[Lemma 8--10]{BB}. We continue to assume $d \asymp T$. 

\begin{lemma}\label{lem3} 
Let $0 < |y| \leq d/30$. Then for any constant $B \geq 0$ one has
\begin{equation}\label{w41}
\Phi_{w_4} (y; d) \ll_{\varepsilon, B} T^{-B}.
\end{equation}
In any case we have   
\begin{equation}\label{w42}
 |y|^j   \Phi^{(j)}_{w_4}(y; d)   \ll_{j, \varepsilon}   T^{3+\varepsilon} (T + |y|^{1/3})^j
 \end{equation}
for any $j \in \Bbb{N}_0$. The same results hold for $\Phi_{w_5}(y; d)$. 
\end{lemma}

The main point of this lemma is the statement that $\Phi_{w_4}(y)$, $\Phi_{w_5}(y)$ decay very quickly, once $|y| \leq T^3$, and that each derivative does not cost more than $T+|y|^{1/3}$, which controls the oscillation. \\

\textbf{Proof.} We insert the   formula \eqref{bessel1} into the definition \eqref{phi} of $\Phi_{w_4}$. We consider the integral over $r$ which is of the form
$$\int_{(0)} F(r) \left( \frac{64\pi^6|y|^2 }{x^3}\right)^r \textbf{spec}^d(r) \frac{dr}{2\pi i}.$$
Integrating by parts, we see that this is negligible unless $x = 4\pi^2|y|^{2/3}(1 + O(T^{-\varepsilon}))$ which we assume from now on. In particular, the argument of the Bessel function in \eqref{bessel1} is $\leq (1+o(1)) \sqrt{8} \pi |y|^{1/3} \leq 10|y|^{1/3}$, and \eqref{w41} follows from \eqref{jbound}. On the other hand, differentiating explicitly with respect to $y$ and estimating trivially using \eqref{jbound} and \eqref{intF}, we obtain \eqref{w42} (since $y/|x| \asymp |y|^{1/3}$ whenever $\Phi_{w_4}$ is non-negligible).\hfill$\square$


 
\begin{lemma}\label{lem4} 
Let $\Upsilon := \min(|y_1|^{1/3} |y_2|^{1/6}, |y_1|^{1/6} |y_2|^{1/3})$. If $\Upsilon \leq d/50$, then
\begin{equation}\label{J1}
\Phi_{w_6}((y_1, y_2); d)   \ll_{B}   T^{-B}
 \end{equation}
for any fixed constant $B \geq 0$. Moreover, 
 \begin{equation}\label{J2}
\begin{split}
|y_1|^{j_1}& |y_2|^{j_2} \frac{\partial^{j_1}}{\partial y_1^{j_1} }\frac{\partial^{j_2}}{\partial y_2^{j_2} } \Phi_{w_6}((y_1, y_2); d) \\
& \ll_{j_1, j_2, \varepsilon } T^{3+\varepsilon}  \left(T+|y_1|^{1/2} + |y_1|^{1/3}|y_2|^{1/6}   \right)^{j_1}  \left(T+ |y_2|^{1/2} + |y_2|^{1/3}|y_1|^{1/6}   \right)^{ j_2} 
\end{split}
 \end{equation}
 for all fixed $j_1, j_2 \in \Bbb{N}_0$. 
\end{lemma}

This lemma says that the cut-off point of $\Phi_{w_6}$ is $  \Upsilon \ll d \asymp T$, and the oscillation is controlled by a factor $T+|y_i|^{1/2} + |y_i|^{1/3}|y_j|^{1/6} $ in the variable $y_i$ (for $j\in \{1, 2\} \setminus \{i\}$).\\

\textbf{Proof.} We insert the Bessel formulae \eqref{bessel} into \eqref{phi} in the various (non-trivial) cases of signs. 
 We consider the integral over of $r$ which is of the form 
 $$\int_{(0)} F(r) \Big| \frac{y_2}{y_1}\Big|^r \left.\begin{cases} x^{3r} (1-x)^{-3r}, & \epsilon= (-, -)\\ x^{-3r}, & \epsilon = (-, +)\\  
x^{3r}, & \epsilon= (+, -) \end{cases} \right\} \textbf{spec}^d(r) \frac{dr}{2\pi i}. $$
Integration by parts shows that this is negligible unless $$x = \begin{cases}  (1 + |y_2/y_1|^{1/3})^{-1} (1 + O(T^{-\varepsilon})), &(\epsilon_1, \epsilon_2) = (-, -),\\ |y_2/y_1|^{1/3}(1 + O(T^{-\varepsilon})), & (\epsilon_1, \epsilon_2) = (-, +), \\    |y_1/y_2|^{1/3}(1 + O(T^{-\varepsilon})), & (\epsilon_1, \epsilon_2) = (+, -), \end{cases}$$
which we assume from now on. In particular, the argument of the first Bessel function in \eqref{bessel} is $$\leq (1 + o(1)) 4\pi \sqrt {|y_1|  + |y_1|^{2/3} |y_2|^{1/3}},$$ while the argument of the second Bessel function   is $$\leq (1+o(1))4\pi\sqrt{ |y_2|  + |y_2|^{2/3} |y_1|^{1/3}}.$$ We can now differentiate $j_1$ times with respect to $y_1$ and $j_2$ times with respect to $j_2$ using \eqref{jder} and estimate trivially using \eqref{jbound} to arrive at \eqref{J2}. The bound \eqref{J1} follows from 
$$4\pi \min \left( \sqrt {|y_1|  + |y_1|^{2/3} |y_2|^{1/3}}, \sqrt{ |y_2|  + |y_2|^{2/3} |y_1|^{1/3}}\right) \leq 4\sqrt{2}\pi \Upsilon \leq 18\Upsilon$$
and \eqref{jbound}. \hfill $\square$ \\

We continue with   bounds  for multiple Fourier transforms. Let $W$ be a smooth weight function with compact support in $(0, \infty)$ all of whose derivatives are uniformly bounded. For $d \geq 2$, $d \in \Bbb{N}$,  $r \in i\Bbb{R}$,  
 $U$,  $V \in \Bbb{R}$, $\Xi \in \Bbb{R} \setminus \{0\}$, we define \begin{equation}\label{deftildeK}
\begin{split}
 &\tilde{\mathcal{K}}_{d, r}(\Xi, U, V):=  \frac{1}{|\Xi|} \int_{0}^{\infty}\int_{0}^{\infty} K_{w_4}( \xi \eta \Xi; d, r) e(\xi U+\eta V) W(\xi) W(\eta) d\xi\, d\eta. 
\end{split}
\end{equation}

\begin{lemma}\label{lem5}  Let   $U, V, \Xi  \in \Bbb{R}$, $|\Xi|, T > 1$.  Let $d \asymp |r| \asymp T$    
and let $\varepsilon, B > 0$.  \\ 
{\rm (a)} We have
$$\tilde{\mathcal{K}}_{d, r}(\Xi, U, 0) \ll_{\varepsilon, B} T^{-B}$$
unless $|U| \leq T^{\varepsilon}$. Similarly, $U = 0$ requires $|V| \leq T^{\varepsilon}$ for a non-negligible contribution.\\
{\rm (b)} We have
$$\tilde{\mathcal{K}}_{d, r}(\Xi, 0, 0) \ll_{\varepsilon, B} \begin{cases} T^{-B}, & |\Xi| \geq T^{3+\varepsilon},\\
T^{-3/2 + \varepsilon}, & \text{always.}\\
\end{cases}$$
{\rm (c)} If $|U|, |V| \geq T^{\varepsilon}$, then $\tilde{\mathcal{K}}_{d, r}(\Xi, U, V) \ll_{\varepsilon, B}  T^{-B}$ unless $U \asymp V$, in which case
$$\tilde{\mathcal{K}}_{d, r}(\Xi, U, V) \ll_{\varepsilon}  (|UV|^{1/2} + T)^{-3/2+\varepsilon}.$$ 
\end{lemma}

 \textbf{Proof.}  We insert the Mellin-Barnes formula \eqref{kw4} getting
\begin{displaymath}
\tilde{\mathcal{K}}_{d, r}(\Xi, U, V):=  \frac{1}{|\Xi|} \int_{0}^{\infty}\int_{0}^{\infty} \int_{(\varepsilon)} |  \xi\eta \Xi |^{1-r-s} G(s; d, r) \frac{ds}{2\pi i}
 e(\xi U+\eta V) W(\xi) W(\eta) d\xi\, d\eta
\end{displaymath} 
with $\varepsilon > 0$ very small and
$$G(s; d, r) = \frac{(\text{sgn}(\Xi)  i)^d}{4\pi^2} |8\pi^3|^{1-r-s} Q(d, s) \Gamma(s+3r) \exp\left(-\frac{\text{sgn}(\Xi) i \pi}{2}(s + 3r)\right),$$
so that 
$$G(s; d, r) \ll   \frac{1}{((d + |s|)(1 + |s + 3r|)^{1/2})^{1-\varepsilon}}, \quad \Re s = \varepsilon.$$
The proof of (a), (b) and (c) is now verbatim the same as in \cite[Lemma 10]{BB}.  \hfill $\square$

\section{Analysis of the Bessel kernels II}

By far the most technical part of the analysis in \cite{BB} was Lemma 11 in that paper, and the present section generalizes this result to the non-spherical case. 

\subsection{Statement of the result} Let $W$ be a smooth weight function with compact support in $(0, \infty)$ all of whose derivatives are uniformly bounded. For $d \geq 2$, $d \in \Bbb{N}$,  $r \in i\Bbb{R}$,  
 $U_1$, $U_2$,   $V_1$, $V_2 \in \Bbb{R}$,   $\Xi_1$, $\Xi_2 \in \Bbb{R} \setminus \{0\}$, we define  \begin{equation}\label{defK}
\begin{split}
\mathcal{K}_{d, r}(\Xi_1, \Xi_2; & U_1, V_1; U_2, V_2):= \frac{1}{|\Xi_1\Xi_2|} \int_{\Bbb{R}^4} K^{\text{sgn}(\Xi_1), \text{sgn}(\Xi_2)}_{w_6}(   \xi_1 \eta_1 \Xi_1,  \xi_2 \eta_2 \Xi_2; d, r) \\
&\times e(\xi_1 U_1+\eta_1 V_1+\xi_2 U_2+\eta_2 V_2)  \   W(\xi_1) W(\eta_1)\overline{W}(\xi_2) \overline{W}(\eta_2)  d\xi_1\, d\xi_2\, d\eta_1\, d\eta_2.\end{split}
\end{equation}

\begin{lemma}\label{lem6} Let  $U_1, V_1, U_2, V_2  \in \Bbb{R}$, $T, |\Xi_1|, |\Xi_2| > 1$, and assume that $\Xi_1, \Xi_2, U_1, U_2, V_1, V_2\ll T^{O(1)}$.   Let $d \asymp |r| \asymp T$  and  
  $\varepsilon, B > 0$. \\
  
\noindent {\rm (a)} We have
 $$\mathcal{K}_{d, r}(\Xi_1, \Xi_2; 0, V_1; U_2, V_2) \ll_{\varepsilon, B} T^{-B}$$
unless $|V_1| \leq T^{\varepsilon}$. Similarly, $V_1 = 0$ requires $|U_1| \leq T^{\varepsilon}$, $U_2 = 0$ requires $|V_2| \leq T^{\varepsilon}$, and $V_2 = 0$ requires $|U_2| \leq T^{\varepsilon}$ for a non-negligible contribution. \\
 
 \noindent {\rm (b)} If $|U_1|, |V_1| \geq T^{\varepsilon}$, we have   
  $$\mathcal{K}_{d, r}(\Xi_1, \Xi_2; U_1, V_1; 0, 0) \ll_{\varepsilon} (|U_1V_1|^{1/4} |\Xi_2|)^{-1+\varepsilon}. $$
  Similarly, if $|U_2|, |V_2| \geq T^{\varepsilon}$,  then  
  $\mathcal{K}_{d, r}(\Xi_1, \Xi_2; 0, 0; U_2, V_2) \ll (|U_2V_2|^{1/4} |\Xi_1|)^{-1+\varepsilon}. $\\

\noindent {\rm (c)} We have $\mathcal{K}_{d, r}(\Xi_1, \Xi_2; 0, 0, 0, 0) \ll_{\varepsilon, B} T^{-B}$ unless $\min(|\Xi_1|, |\Xi_2|) \geq T^{3-\varepsilon}$. In this case we have
\begin{equation}\label{c1}
\mathcal{K}_{d, r}(  \Xi_1,  \Xi_2; 0, 0; 0, 0) \ll_{\varepsilon} T^{-3+\varepsilon}
\end{equation}
and 
\begin{equation}\label{c2}
\sum_{\epsilon \in \{\pm 1\}^2} \sum_D \frac{\phi(D)}{D^2} \mathcal{K}_{d, r}\left(\frac{\epsilon_1 \Xi_1}{D}, \frac{\epsilon_2 \Xi_2}{D}; 0, 0; 0, 0\right) \ll_{\varepsilon} |\Xi_1 \Xi_2|^{-1/2 + \varepsilon}.
\end{equation}

\noindent {\rm (d)} For $|U_1|, |V_1|, |U_2|, |V_2| \geq T^{\varepsilon}$ we define
\begin{equation}\label{defw}
\Upsilon_1 = \frac{|\Xi_1|}{   |U_1V_1|}, \quad \Upsilon_2 = \frac{|\Xi_2|}{    |U_2V_2|}.
\end{equation}
Let $G$ be a  smooth function with compact support all of whose derivatives are uniformly bounded and let $\rho_0\in \Bbb{R}$ with $|\rho_0| \asymp T$. 
 Then we 
have  
\begin{equation}\label{11cbound1}
|U_1V_1U_2V_2|^{1/2}  \int_{\Bbb{R}} G\left(\frac{\rho-\rho_0}{T^{\varepsilon}}\right)  \mathcal{K}_{d, i\rho}(\Xi_1, \Xi_2; U_1, V_1; U_2, V_2) 
d\rho \ll_{\varepsilon}     T^{-1 +\varepsilon}
\end{equation}
 if
\begin{equation}\label{specialcase}
 |\Upsilon_1 - 1| + \Upsilon_2 \ll T^{-\Halfbfrac} \ll   |U_2V_2|T^{-2}    \quad    \text{or} \quad  |\Upsilon_2 - 1| + \Upsilon_1 \ll T^{-\Halfbfrac} \ll  |U_1V_1|  T^{-2},
 \end{equation}
and 
\begin{equation}\label{11cbound2}
|U_1V_1U_2V_2|^{1/2} \int_{\Bbb{R}} G\left(\frac{\rho-\rho_0}{T^{\varepsilon}}\right) \mathcal{K}_{d, i\rho}(\Xi_1, \Xi_2; U_1, V_1; U_2, V_2)
d\rho  \ll_{\varepsilon}    T^{-1 -\Eighthbfrac + \varepsilon}
\end{equation}
otherwise. \\

\noindent {\rm (e)} Keep the assumptions and notation from part (d) and assume   in addition
\begin{equation}\label{condition-e}
|U_2V_2| \geq (|U_1V_1| + T^2)T^{\varepsilon},
\end{equation}
then  $$\int_{\Bbb{R}} G\left(\frac{\rho-\rho_0}{T^{\varepsilon}}\right)  \mathcal{K}_{d,  i\rho}(\Xi_1, \Xi_2; U_1, V_1; U_2, V_2)
d\rho
   \ll_{\varepsilon, B} T^{-B}$$ unless 
\begin{equation}\label{cond11c}
|\Upsilon_2 - 1| \leq \frac{T^{\varepsilon} }{|U_2V_2|^{1/4}} + \frac{|U_1 V_1|^{1/2}+T}{|U_2V_2|^{1/2}}.
\end{equation}
Similarly,  $|U_1V_1| \geq (|U_2V_2| + T^2)T^{\varepsilon}$ requires 
$|\Upsilon_1 - 1| \leq  T^{\varepsilon}|U_1V_1|^{-1/4} +(|U_2 V_2|^{1/2}+T)|U_1V_1|^{-1/2}$
for a non-negligible contribution. 
\end{lemma} 

\emph{Remark:} This is almost literally the analogue of \cite[Lemma 11]{BB} with one important exception. In parts (d) and (e) we have inserted an additional short integration over $\rho$ of length roughly $T^{\varepsilon}$. Although  it is primarily the purpose of Lemma \ref{lem4} to exploit the integration over $\rho = \Im r$ coming from the definition \eqref{phi}, remembering the $\rho$-integration  also in the present lemma will produce one extra condition from partial integration that simplifies the proof   compared to \cite[Lemma 11]{BB} (although it still remains lengthy and complicated, and for technical reasons the numerical values of the exponents change). 

\subsection{Proof of parts (a), (b), (c)}
 Again we insert the Mellin-Barnes representation \eqref{kw6}, and for convenience we change variables $s_1 \mapsto s_1 - r$, $s_2 \mapsto s_2 + r$.  
 We must be a little careful with the choice of the $s_1, s_2$-contours, but as in \cite[Lemma 11]{BB} we see that we can integrate over $s_1, s_2 \in [\varepsilon - iT^B, \varepsilon + iT^B]$ for some sufficiently large constant $B$, up to a negligible error. Hence $4\pi^2  \mathcal{K}_{d, r}(\Xi_1, \Xi_2; U_1, V_1; U_2, V_2)$ equals
 \begin{equation}\label{analyzing}
\begin{split}
 \frac{1}{|\Xi_1\Xi_2|} \int_{\Bbb{R}^4} & \int_{|t_j| \leq T^B}  | 4\pi^2  \xi_1\eta_1\Xi_1|^{1-\varepsilon - it_1} |   4\pi^2 \xi_2\eta_2\Xi_2|^{1-\varepsilon - it_2}  G^{\epsilon}((s_1, s_2); d,  r)  \frac{dt_1\, dt_2 }{(2\pi)^2} \\
&\times e(\xi_1 U_1+\eta_1 V_1+\xi_2 U_2+\eta_2 V_2)  \   W(\xi_1) W(\eta_1)\overline{W}(\xi_2) \overline{W}(\eta_2)  d\xi_1\, d\xi_2\, d\eta_1\, d\eta_2
\end{split}
\end{equation}
(up to a negligible error) with
\begin{equation}\label{g}
G^{\epsilon}((s_1, s_2); d,  r) =     B_{w_6}^{\epsilon}((s_1 - r, s_2+r ), r) Q(d, s_1-r) Q(d, s_2+r)  
\end{equation}
where $\epsilon = (\text{sgn}(\Xi_1), \text{sgn}(\Xi_2)).$ We recall from \eqref{B} that
$$B_{w_6}^{\epsilon}((s_1-r, s_2+r), r) = \begin{cases} (-1)^d B(s_1 + 2r, s_2 - 2r), 
\\ B(s_2-2r, 1 - s_1 - s_2), 
\\ B(s_1+2r, 1 - s_1 - s_2)%
\end{cases}$$
in the various non-trivial cases of signs. The function
$G^{\epsilon}((s_1, s_2); d, r)$ now plays exactly the role of $G((s_1, s_2), \mu)S^{\epsilon_1, \epsilon_2}((s_1, s_2), \mu)$ in the proof of \cite[Lemma 11]{BB}. \\

The proofs of parts (a) -- (c) are now verbatim the same as in \cite[Section 15.1]{BB}; for (b) we use that\footnote{Note the typographical difference between $\epsilon \in \{\pm 1\}$ or $\{\pm 1\}^2$ and $\varepsilon > 0$ here and throughout the paper.} 
$$\int_{t_1 \asymp U_1} \int_{|t_2| \leq T^{\varepsilon}} |G^{\epsilon}((\sigma_1 + it_1, \sigma_2 + it_2); d,  r)|   dt_1\, dt_2 \ll T^{\varepsilon} U_1  \frac{(T+|U_1|)^{3\sigma_1 - 3/2} T^{3\sigma_2 - 3/2}}{|U_1|^{\sigma_1 + \sigma_2 - 1/2}}$$
for $0 < \sigma_1 < 1/2 < \sigma_2 < 1$ by Stirling's formula. Since (c) is rather delicate, we give some more details. For $U_1 = V_1 = U_2 = V_2 = 0$ we have  
$$4 \pi^2 \mathcal{K}_{d, r}(\Xi_1, \Xi_2; 0, 0; 0, 0) =   \int_{-i\infty}^{i\infty}\int_{-i\infty}^{i\infty}  G^{\epsilon} ((s_1, s_2); d, r)  |  \Xi_1|^{-s_1}|  \Xi_2|^{-s_2} \widehat{ {\tt W}}(2-s_1) \widehat{\overline{ {\tt W}}}(2-s_2)  \frac{ds_1 \,  ds_2}{(2\pi i)^2},$$
where $${\tt W}(x) = \frac{1}{4\pi^2}  \int_{\Bbb{R}} W(\xi) W\left(\frac{x}{4\pi^2\xi}\right) \frac{d\xi}{\xi}.$$
Here we may shift the $s_1$ and/or $s_2$ contour to the left. The poles at $s_1 = - n -2r$, $s_2 = - n +2r$ for $n \in \Bbb{N}_0$ contribute negligibly by the rapid decay of $\widehat{ {\tt W}}(s)$ and $|r| \asymp T$. By Stirling's formula we see that this forces $|\Xi_1|$ and $|\Xi_2|$ to be $\gg T^{3-\varepsilon}$ for a non-negligible contribution, and we obtain \eqref{c1} for the contours at $\Re s_1 = \Re s_2 = \varepsilon$. In order to prove \eqref{c2}, we shift both contours to $\Re s_1 = \Re s_2 = -1/2$   getting
\begin{displaymath}
\begin{split}
&4\pi^2 \sum_{\epsilon \in \{\pm 1\}^2}   \sum_D \frac{\phi(D)}{D^2} \mathcal{K}_{d, r}\left(\frac{\epsilon_1\Xi_1}{D}, \frac{\epsilon_2\Xi_2}{D}; 0, 0; 0, 0\right) \\
&= \sum_{\epsilon \in \{\pm 1\}^2}   \int_{(-1/2)} \int_{(-1/2)} \frac{\zeta(1 - s_1 - s_2)}{\zeta(2 - s_1 - s_2)} G^{\epsilon}((s_1, s_2); d, r) |  \Xi_1|^{-s_1}|  \Xi_2|^{-s_2} \widehat{ {\tt W}}(2-s_1) \widehat{\overline{ {\tt W}}}(2-s_2)  \frac{ds_1\, ds_2}{(2\pi i)^2}
\end{split}
\end{displaymath}
plus a negligible error. The key observation is that 
$$\sum_{\epsilon \in \{\pm 1\}^2} G^{\epsilon}((s, -s); d, r) = 0.$$
Indeed, we have  
$$B_{w_6}^{(-, -)}((s-r, -s+r), r)  = B_{w_6}^{(+, -)}((s-r, -s+r), r) + B_{w_6}^{(-, +)}((s-r, -s+r), r) = 0.$$
Therefore we may shift to the right past the removable pole at $s_1 + s_2 = 0$  to $\Re s_1 = \Re s_2 = 1/2 - \varepsilon$ and conclude \eqref{c2} from Stirling's formula. 

We note in passing that as in \cite[Lemma 11]{BB}, although in a slightly different form, a combination of signs is necessary to create a zero that is absolutely crucial for the success of the proof. 

\subsection{The general set-up for (d) and (e)} The rest of this section is devoted to the proof of parts (d) and (e), where we assume $|U_1|, |U_2|, |V_1|, |V_2| \geq T^{\varepsilon}$.    We return to \eqref{analyzing} and integrate over $\xi_1, \xi_2, \eta_1, \eta_2$ using the stationary phase lemma \cite[Lemma 3]{BB}. In particular, we see that this 4-fold integral is negligible unless
$$t_1 \asymp U_1 \asymp V_1, \quad t_2 \asymp U_2\asymp V_2$$
in which case    $\mathcal{K}_{d, r}(\Xi_1, \Xi_2; U_1, V_1; U_2, V_2)$ equals (up to a negligible error)
\begin{equation}\label{k-mb}
\frac{1}{|U_1V_1U_2V_2|^{1/2}} \int_{\Bbb{R}^2}  G^{\epsilon}((\varepsilon+ it_1, \varepsilon+it_2); d,  r) (e^{-2}\Upsilon_1t_1^2)^{-it_1}(e^{-2}\Upsilon_2t_2^2)^{-it_2} \Psi\left(\frac{t_1}{U_1}, \frac{t_2}{U_2}\right) dt_1\, dt_2
\end{equation}
where $\Upsilon_j$ is as in \eqref{defw} 
and $\Psi$ is a  smooth, compactly supported function all of whose derivatives are uniformly bounded. For notational simplicity we write $\rho = \Im r$. 

Stirling's formula is easy to apply for the $Q$-factors in \eqref{g} since $d$ is large, but we need to insert a partition of unity to treat the beta function if $\Im s_1$ or $\Im s_2$ are close to $-2r$ resp.\ $2r$, or if $s_1+s_2$ is  close to 0. Thus we insert a localizing factor
$$F_1\left(\frac{t_1 + 2\rho}{B_1}\right) F_2\left(\frac{t_2 - 2\rho}{B_2}\right)  F_3\left(\frac{t_1 + t_2}{B_3}\right)$$
where $B_1, B_2, B_3 \in \{2^\nu \mid \nu \in \Bbb{N}_0\}$ and $F_j$ are fixed smooth functions with compact support in $[-2, 2]$ if the corresponding $B_j = 1$ and support in $[1/2, 2] \cup [-2, -1/2]$ if $B_j \geq 2$. The triangle inequality implies that we can have only 
\begin{equation}\label{BB}
B_1 \asymp B_2 \geq B_3 \quad \text{or} \quad B_1 \leq B_2 \asymp B_3 \quad \text{or} \quad B_2 \leq B_1 \asymp B_3,
\end{equation}
the last two cases being symmetric. 
Stirling's formula implies that the integrand in \eqref{k-mb} is
$$\ll T^{\varepsilon} \frac{1}{ (T+ |U_1|) (T+ |U_2|)}  \left(\frac{ B_3}{B_1B_2}\right)^{1/2}.$$
  
If two of the three values $B_1, B_2, B_3$ are $\ll T^{\varepsilon}$, then all three satisfy this bound, and a trivial estimate shows that  $\mathcal{K}_{d, r}(\Xi_1, \Xi_2; U_1, V_1; U_2, V_2) \ll T^{-2+\varepsilon}|U_1V_1U_2V_3|^{-1/2}$, which easily implies a stronger version of \eqref{11cbound1} and \eqref{11cbound2}. Moreover,   the condition   \eqref{condition-e} is void. If only $B_3 \leq T^{\varepsilon}$, then again a trivial estimate shows $\mathcal{K}_{d, r}(\Xi_1, \Xi_2; U_1, V_1; U_2, V_2) \ll T^{-3/2+\varepsilon}|U_1V_1U_2V_3|^{-1/2}$ which still  implies a stronger version of \eqref{11cbound1} and \eqref{11cbound2}, and  the  condition \eqref{condition-e} is void. In particular, for the proof of parts (d) and (e) we can assume that  (say) $ B_2, B_3 \geq T^{\varepsilon}$. For now we will also assume $B_1 \geq T^{\varepsilon}$ and treat the case $B_1 \leq T^{\varepsilon}$ on the way. (Note that under the symmetric assumption $B_2 \leq T^{\varepsilon}$  the condition \eqref{condition-e} would be void, so that the present assumptions are really no loss of generality.)

We now insert Stirling's formula for $G^{\epsilon}((\varepsilon+ it_1, \varepsilon+it_2); d,  r)$. As in \cite[(15.4)]{BB}, we can express $|U_1V_1U_2V_2|^{1/2}\mathcal{K}_{d, r}(\Xi_1, \Xi_2; U_1, V_1; U_2, V_2)$ up to a negligible error as sums over integrals of the form
$$\mathcal{I}(B_1, B_2, B_3) = \int_{\Bbb{R}^2} e^{ig(t_1, t_2; d, r)} \mathcal{F}(t_1, t_2; d, r) dt_1\, dt_2$$
for various choices of $B_1, B_2, B_3$, where a standard application of Stirling's formula gives
\begin{displaymath}
\begin{split}
g(t_1, t_2; d, r)= &d \arctan\left(\frac{t_1 - \rho}{d/2}\right) +  (t_1 - \rho) \log \frac{ (t_1 - \rho)^2 + (d/2)^2}{e^2}\\
& + d \arctan\left(\frac{t_2 + \rho}{d/2}\right) +   (t_2 + \rho) \log \frac{ (t_2 + \rho)^2 + (d/2)^2}{e^2}\\
&+  (t_1 + 2\rho) \log \frac{|t_1 + 2\rho| }{e} +  (t_2 - 2\rho) \log \frac{|t_2 - 2\rho| }{e} -  (t_1 + t_2) \log \frac{|t_1 + t_2| }{e}\\
&- t_1\log t_1^2 e^{-2}\Upsilon_1 - t_2 \log t_2^2 e^{-2}\Upsilon_2
\end{split}
\end{displaymath}
and 
 $$\frac{\partial^k}{\partial r^k} \frac{\partial^n}{\partial t_1^n} \frac{\partial^m}{\partial t_2^m}  \mathcal{F}(t_1, t_2; d,  r) \ll_{n, m, k} \frac{T^{\varepsilon} }{ (T+ |U_1|) (T+ |U_2|)}  \left(\frac{B_3}{B_1B_2}\right)^{1/2} \frac{1}{E_1^nE_2^mF^k} $$ 
for $n, m, k \in \Bbb{N}_0$ with
$$E_1 := \min(B_1, B_3, |U_i|), \quad  E_2 := \min(B_2, B_3, |U_i|), \quad F  := \min(B_1, B_2, T).$$

It is convenient to introduce the notation
$$A \lll B$$
to mean $A \leq \delta B$ for a sufficiently small  constant $\delta$ 
  (where ``sufficiently small'' depends on the implicit constants in the condition $d \asymp |r| \asymp T$, the support of the weight functions and $\varepsilon$ where applicable). Similarly we write $A \ggg B$ to mean $A \geq \Delta B$ for a sufficiently large constant $\Delta$.

\subsection{Computations with derivatives} 
We compute
\begin{align}
  \nonumber &g_1(t_1, t_2; d, r) :=  \frac{\partial}{\partial t_1} g(t_1, t_2; d, r)  = \log\Big| \frac{(t_1 + 2\rho)((d/2)^2 + (\rho-t_1)^2)}{ (t_1+t_2)t_1^2 \Upsilon_1}\Big| ,\\  \label{7.14}
&g_2(t_1, t_2; d, r) :=  \frac{\partial}{\partial t_2} g(t_1, t_2; d, r)   = \log\Big| \frac{(t_2 - 2\rho)((d/2)^2 + (\rho+t_2)^2)}{ (t_1+t_2)t_2^2 \Upsilon_2} \Big|,\\   \label{7.15}
& h(t_1, t_2; d, r) :=  \frac{\partial}{\partial r} g(t_1, t_2; d, \rho)   = \log\Big| \frac{((d/2)^2 + (\rho+t_2)^2)(t_1 + 2\rho)^2 }{((d/2)^2 + (\rho-t_1)^2)(t_2 - 2\rho)}  \Big|.
\end{align}
It is easy to see that
\begin{equation}\label{156}
\frac{\partial}{\partial t_2} g_1(t_1, t_2; d, r) \ll \frac{1}{B_3}, \quad \frac{\partial}{\partial t_1} g_2(t_1, t_2; d, r) \ll \frac{1}{B_3}
\end{equation}
and
\begin{equation}\label{156b}
\frac{\partial^n}{\partial t_i^n} g_i(t_1, t_2; d, r) \ll_n \frac{1}{E_i}\end{equation}
for $i = 1, 2$ and $n \in \Bbb{N}$. 
We also have
\begin{equation}\label{alsohave}
\frac{\partial}{\partial t_1} g_1(t_1, t_2;d, r) = - \frac{2}{t_1} + \frac{1}{t_1 + 2\rho} - \frac{1}{t_1 + t_2} + \frac{2(t_1 - \rho)}{(d/2)^2 + (t_1 - \rho)^2} = \frac{t_2t_1^3 + \ldots}{t_1^5 + \ldots},
\end{equation}
where $\ldots$ in the numerator are homogeneous expressions in $t_1, t_2, d, \rho$ of total degree 4 whose highest power of $t_1$ is 2, and $\ldots$ in the denominator are homogeneous expressions in $t_1, t_2, d, \rho$ of total degree 5 whose highest power of $t_1$ is 4.  In particular, if $|U_1| \ggg T + |U_2|$, then inductively we see that 
\begin{equation}\label{stronger}
\frac{\partial^n}{\partial t_1^n} g_1(t_1, t_2; d, r) = \frac{(-1)^{n-1} n! t_2t_1^{4n-1} + \ldots}{t_1^{5n} + \ldots} \asymp \frac{|U_2|}{|U_1|^{n+1}}
\end{equation}
for each fixed $n \in \Bbb{N}$. A similar relation holds with exchanged indices in the case $|U_2| \ggg T + |U_1|$.  Note that the present version of \cite[(15.13) \& (15.15)]{BB} is a little stronger since the condition $|U_1U_2| \ggg T^2$ is not needed. 

Finally, if $B_2 \lll B_1$, we obtain from the exact formula in  \eqref{alsohave} that 
\begin{equation}\label{B2small}
 \frac{\partial}{\partial t_1} g_1(t_1, t_2; d, r) =  - \frac{2((d/2)^2 + \rho^2  - \rho t_1)}{t_1((d/2)^2 + (\rho- t_1)^2)} + O\left(\frac{B_2}{B_1^2}\right)
\end{equation}
and similarly 
\begin{equation}\label{B1small}
\frac{\partial}{\partial t_2} g_2(t_1, t_2; d, r)  = - \frac{2((d/2)^2 + \rho^2 + \rho t_2)}{t_2((d/2)^2 + (\rho+t_2)^2)} + O\left(\frac{B_1}{B_2^2}\right)
\end{equation}
if $B_1 \lll B_2$. We emphasize that \eqref{B1small} is applicable even in the previously excluded case (which is precisely where we need it) when $B_1 \leq T^{\varepsilon}$, where we cannot insert Stirling's formula for $\Gamma(s_1 + 2r)$, but this factor is irrelevant in the computation of $(\partial/\partial t_2) g_2$. Similarly, \eqref{B2small} is applicable in the symmetric case $B_2 \leq T^{\varepsilon}$. 

Integrating by parts with respect to $t_1, t_2$ using \eqref{156b}, we see as in \cite[(15.16)]{BB} that $\mathcal{I}(B_1, B_2, B_3)$ is negligible  unless 
\begin{equation}\label{e1/2}
g_i(t_1, t_2; d, r) \ll T^{\varepsilon} E_i^{-1/2}.
\end{equation}
In special situations $|U_1| \ggg T + |U_2|$ resp.\ $|U_2| \ggg T + |U_1|$ as considered in \eqref{stronger} and the subsequent paragraph we have the slightly stronger bounds
\begin{equation}\label{u1u2}
g_1(t_1, t_2; d, r) \ll T^{\varepsilon} |U_2|^{1/2} |U_1|^{-1} \quad \text{resp.} \quad g_2(t_1, t_2; d, r) \ll T^{\varepsilon} |U_1|^{1/2} |U_2|^{-1}.
\end{equation}

 In particular, we can assume the consistency relations
\begin{equation}\label{consist}
\frac{B_1(T+|U_1|)^2}{ B_3|U_1|^2 \Upsilon_1} \asymp 1 \asymp  \frac{ B_2(T+|U_2|)^2}{  B_3|U_2|^2 \Upsilon_2}.
\end{equation}
This is \cite[(15.18)]{BB}. However,  statements (d) and (e) of the lemma also involve an integral over $\rho$ of length $T^{\varepsilon}$, and partial integration in $\rho$ along with \eqref{7.15} gives the additional consistency relation
\begin{equation}\label{r}
1 \asymp \frac{(T^2 + U_2^2)B_1^2}{(T^2 + U_1^2) B_2^2}.
\end{equation}
This is the only place where the extra integration over $\rho$ is used. Once we have recorded this condition, we will not use the $\rho$-integration any more and treat $\rho$ as a fixed number satisfying $|\rho| \asymp T$. 

We now state a slightly simplified version of \cite[Sublemma 1]{BB}. Let $\mathcal{A}$ be the set of all $(t_1, t_2)$ such that
$$   t_i \asymp U_i, \quad |t_1 + 2\rho| \asymp B_1, \quad |t_2 -2 \rho| \asymp B_2, \quad |t_1+t_2| \asymp B_3$$
and $g_i(t_1, t_2; d, r)$ satisfies \eqref{e1/2} or \eqref{u1u2} if the respective conditions are satisfied. We have
\begin{equation}\label{basicbound}
\mathcal{I}(B_1, B_2, B_3) \ll T^{\varepsilon}  \frac{1}{ (T+ |U_1|) (T+ |U_2|)}  \left(\frac{ B_3}{B_1B_2}\right)^{1/2} \text{meas}(\mathcal{A}).
\end{equation}
\begin{lemma}\label{sublemma} Let $M_1,M_2 \subseteq \Bbb{R}$ be two, possibly infinite, intervals. 
Define $M = M_1 \times M_2$,
\begin{align*}
	H_1 = \inf\left\{\left|\frac{\partial g_1}{\partial t_1} (t_1, t_2; d, r)\right| : (t_1, t_2) \in \mathcal{A}  \cap M \right\},  \quad H_2 = \inf\left\{\left|\frac{\partial g_2}{\partial t_2} (t_1, t_2; d, r)\right| : (t_1, t_2) \in \mathcal{A}  \cap M\right\},
\end{align*}
We have the following estimates for the measure of $\mathcal{A}$:
\begin{equation}\label{1a}
 \text{{\rm meas}}(\mathcal{A}  \cap M) \ll  T^{\varepsilon} \frac{\min(\text{{\rm meas}}(M_1), B_1, |U_1|)}{\sqrt{E_2} H_2}, \quad \text{provided } H_2 \gg \frac{T^{\varepsilon}}{E_2^{5/4}};
\end{equation}
 \begin{equation}\label{1b}
 \text{{\rm meas}}(\mathcal{A} )  \ll T^{\varepsilon} \frac{\min(B_1, |U_1|)B_3}{\sqrt{E_1} }; 
\end{equation}
\begin{equation}\label{2a}
 \text{{\rm meas}}(\mathcal{A} \cap M) \ll T^{\varepsilon} \frac{\min(\text{{\rm meas}}(M_2), B_2, |U_2|)}{\sqrt{E_1} H_1}, \quad \text{provided } H_1 \gg \frac{T^{\varepsilon}}{E_1^{5/4}};
\end{equation}
 \begin{equation}\label{2b}
 \text{{\rm meas}}(\mathcal{A}  ) \ll  T^{\varepsilon} \frac{\min(  B_2, |U_2|)B_3}{\sqrt{E_2} }. 
\end{equation}
 If in addition $|U_1| \ggg T + |U_2|$, then we have
\begin{equation}\label{modif1}
\text{{\rm meas}}(\mathcal{A}) \ll |U_2|^{-1/2  } B_2 |U_1|^{1+\varepsilon};
\end{equation}
similarly, if $|U_2| \ggg T + |U_1|$, then we have
\begin{equation}\label{modif2}
 \text{{\rm meas}}(\mathcal{A}) \ll |U_1|^{-1/2  } B_1 |U_2|^{1+\varepsilon}.
\end{equation}
\end{lemma}
We emphasize that the more complicated bound \cite[(15.23)]{BB} whose proof requires a two-dimensional Taylor argument,  is not needed, because we have the extra relation \eqref{r} instead.

The lemma is proved exactly as in \cite[Sublemma 1]{BB} and we refer to this paper for   more details.   By a bit of Morse theory, the number of connected components of $\mathcal{A}$ is uniformly bounded, so it suffices to bound the measure of each connected component.   The rest  of that proof is an   application of Taylor's theorem and the above bounds \eqref{156} for the derivatives. 

\subsection{Proof of part (e)} 
If $|U_2| \gg (|U_1| + T) T^{\varepsilon}$, then by \eqref{7.14} we have
$$g_2(t_1, t_2; d, r) =   - \log \Upsilon_2 + O\left(\frac{|U_1| + T}{|U_2|}\right).$$
Since   $B_2, B_3 \asymp |U_2|$, we have $E_2 \asymp |U_2|$, so that in view of \eqref{e1/2} we obtain
$$|\Upsilon_2 - 1| \ll  \frac{|U_1| + T}{|U_2|} + \frac{T^{\varepsilon}}{|U_2|^{1/2}}, $$
which in view of $|U_2| \asymp |V_2|$ is equivalent   to  \eqref{cond11c}.  Notice that for this argument we do not need to insert Stirling's formula for $\Gamma(s_1 +2r)$, so that the argument works even in the previously excluded case $B_1 \leq T^{\varepsilon} \leq B_2, B_3$, hence the proof of (e) is complete.

\subsection{The case where $B_1$ is small} We are now prepared to treat the remaining exceptional case where (say) $B_1 \leq T^{\varepsilon}$. This implies in particular $|U_1| \asymp T$ and $B_2 \asymp B_3$. We distinguish several cases depending on the size of $|U_2|$. 

 \emph{Case 1:} Suppose that $|U_2| \ggg T$, so that $B_2 \asymp B_3 \asymp E_2 \asymp |U_2|$. In this case we apply \eqref{modif2} with $B_1 \leq T^{\varepsilon}$ and  \eqref{basicbound} to  obtain
$$\mathcal{I}(B_1, B_2, B_3) \ll \frac{ T^{\varepsilon} }{ (T+ |U_1|) (T+ |U_2|)}  \left(\frac{B_3}{B_2}\right)^{1/2} \frac{|U_2|}{T^{1/2}} \ll T^{-3/2 + \varepsilon}.$$

\emph{Case 2:} Suppose that $|U_2| \asymp  T$. Then $B_2 \asymp B_3 \ll T$. If $B_2 \leq T^{9/10}$, we can estimate trivially $\text{meas}(\mathcal{A}) \ll B_1B_2 \leq T^{\varepsilon} B_2$, so that 
$$\mathcal{I}(B_1, B_2, B_3)  \ll T^{\varepsilon} \frac{1}{T^2} \cdot B_2 \leq T^{-1 - \frac{1}{10} + \varepsilon}.$$
Suppose from now on that $B_2   \geq T^{9/10}$. Formula \eqref{B1small} tells us that typically $(\partial/\partial t_2)g_2 \asymp  1/|U_2|$, so that we can apply \eqref{case1a}. Unfortunately  there may be a small region of $t_2$ where the numerator in \eqref{B1small} could drop.  Therefore in this case, we let $M_2$ be the region
$$\Big|t_2 + \frac{(d/2)^2 + \rho^2}{\rho}\Big| \geq  T^{9/10}, $$
so that $\mathcal{I}(B_1, B_2, B_3)$ is a sum of integrals $\mathcal{I}_{M_2}$ and $\mathcal{I}_{\Bbb{R}\setminus M_2}$ over $M_2$ and its complement, respectively. For the contribution of $M_2$, we have $(\partial/\partial t_2)g_2 \gg  T^{-1/10}|U_2|^{-1} \gg E_2^{-5/4+\varepsilon}$, so that \eqref{case1a} is applicable and gives $\text{meas}(\mathcal{A} \cap (\Bbb{R} \times M_2)) \ll T^{11/10+\varepsilon} E_2^{-1/2}$. We obtain  $$\mathcal{I}_{M_2} \ll T^{\varepsilon}\frac{1}{T^2} \frac{1}{T^{1/2 \cdot 9/10} T^{-11/10}} = T^{-27/20 + \varepsilon}.$$
The contribution of the complement can be estimated trivially to be $$\mathcal{I}_{\Bbb{R}\setminus M_2} \ll T^{\varepsilon} \frac{1}{T^2}  T^{9/10} = T^{-11/10  +\varepsilon}.$$

\emph{Case 3:} Finally suppose $|U_2| \lll T$. Here we always have $(\partial/\partial t_2)g_2 \asymp  1/|U_2|$ as well as $B_2 \asymp B_3 \asymp T$, $E_2 \asymp |U_2|$. By \eqref{1a} we obtain easily
$$\mathcal{I}(B_1, B_2, B_3) \ll T^{\varepsilon - 2} \text{meas}(\mathcal{A})\ll T^{\varepsilon - 2} |U_2|^{1/2} \ll  T^{-3/2 + \varepsilon} .$$

This proves a numerically stronger version of \eqref{11cbound1} and \eqref{11cbound2} in all  cases and completes the discussion of the case $B_1 \leq T^{\varepsilon}$.

\subsection{The nearly generic case}
We choose some small positive constants $0 < u_1, u_2, b, c, d < 1/5$ satisfying 
\begin{equation}\label{exponents}
u_1 > b, \quad u_2   > u_1 + 12b + d, \quad  c > 6b, \quad d > 8b, 
\end{equation}
 and we consider in this subsection  the case
\begin{equation}\label{case1a}
T^{1-b} \leq |U_1|, |U_2|, B_1, B_2, B_3 \leq T^{1+b}.
\end{equation}
 Our assumption implies $T^{1-b} \leq E_1, E_2 \leq  T^{1+b}$, and from the consistency relation \eqref{consist} we also have $T^{-b} \ll \Upsilon_l \ll T^{3b}$, $l=1,2$. As in \cite[Section 15.6]{BB}, we want to analyze $\mathcal{A}$ more directly. The conditions $g_i(t_1, t_2; d, r) \ll T^{\varepsilon} E_i^{-1/2}$ are equivalent to 
 \begin{equation}\label{cubic}
\begin{split}
&(\Upsilon_1^{-1} - \alpha_1)t_1^3 - \Upsilon_1^{-1}C_1 t_1 - \Upsilon_1^{-1} C_2 - \alpha_1 t_1^2t_2 \ll T^{\varepsilon} B_3 |U_1|^2 E_1^{-1/2},\\
& (\Upsilon_2^{-1} - \alpha_2)t_2^3 - \Upsilon_2^{-1}C_1 t_2 + \Upsilon_2^{-1} C_2 - \alpha_2 t_2^2t_1 \ll T^{\varepsilon} B_3 |U_2|^2 E_2^{-1/2}.
\end{split}
\end{equation}
for $\alpha_1, \alpha_2 \in \{\pm 1\}$ and 
$$C_1 =  3\rho^2 - (d/2)^2 \ll T^2, \quad C_2 = -2\rho((d/2)^2 + \rho^2) \asymp \pm T^3.$$
It is at this point where things become more complicated than in \cite{BB}, since the size of $C_1$ can drop and be much smaller than its generic size $T^2$. To begin with, we follow the argument in \cite[(15.29)]{BB}, solve the first condition for $t_2$ up to an error of size $O(T^{\varepsilon} B_3E_1^{-1/2})$ and substitute this into the second, getting
\begin{align*}
	\sum_{i=0}^9 a_i t_1^i  = a_9\prod_{i=1}^9 (t_1 - q_i)  \ll& T^{\varepsilon} B_3 |U_1|^6 \left(E_1^{-1/2}((1 + \Upsilon_2^{-1})|U_2|^2 + \Upsilon_2^{-1} T^2 + |U_1U_2|) + E_2^{-1/2}|U_2|^2\right) \\
	\ll & T^{\frac{17}{2}+\frac{21}{2}b+\varepsilon}
\end{align*}
for some complex numbers $a_i$, $q_i$, independent of $t$, where in particular
\begin{displaymath}
\begin{split}
a_9 &= (\Upsilon_1^{-1} - \alpha_1)^2 (\Upsilon_1^{-1}\Upsilon_2^{-1} -\alpha_2\Upsilon_1^{-1} - \alpha_1\Upsilon_2^{-1}),\\
a_8 & = 0,\\
a_7 & = -C_1(\Upsilon_1^{-1} - \alpha_1) (\alpha_1\alpha_2 \Upsilon_1^{-1} - 3\alpha_2 \Upsilon_1^{-2} + \Upsilon_2^{-1} - 3\alpha_1 \Upsilon_1^{-1} \Upsilon_2^{-1} + 3\Upsilon_1^{-2} \Upsilon_2^{-1}),\\
a_6 & = - C_2(-\alpha_2\Upsilon_1^{-1} + 4\alpha_1\alpha_2 \Upsilon_1^{-2} - 3\alpha_2 \Upsilon_1^{-3} - \alpha_1 \Upsilon_2^{-1} + 3\Upsilon_1^{-1} \Upsilon_2^{-1} - 6\alpha_1 \Upsilon_1^{-2} \Upsilon_2^{-1} + 3\Upsilon_1^{-3} \Upsilon_2^{-1}).
\end{split}
\end{displaymath}

\emph{Case I:}  Assume that $|\Upsilon_1^{-1} - \alpha_1| \geq T^{-u_1}$ and $|\Upsilon_1^{-1}\Upsilon_2^{-1} - \alpha_2\Upsilon_1^{-1} - \alpha_1 \Upsilon_2^{-1}| \geq T^{-u_2}$. 
We conclude that
$$|t_1 - q_i|^9 \ll \frac{T^{\frac{17}{2}+\frac{21}{2}b+\varepsilon}}{|a_9|} \leq T^{\frac{17}{2}+2u_1 + u_2+\frac{21}{2}b+\varepsilon}, $$
for some $i \in \{1, \ldots, 9\}.$
Since $t_1$ is now in a fixed interval, independent of $t_2$, we have also determined $B_2$ within an interval (depending on $t_2$) of length $O(T^{\varepsilon} B_3E_1^{-1/2})$, so that 
\begin{equation}\label{lem2-1}
\text{{\rm meas}}(\mathcal{A}  ) \ll T^{\frac{17}{18}+\frac{2u_1}{9} + \frac{u_2}{9}+\frac{7}{6}b+\varepsilon} B_3 E_1^{-1/2} \leq T^{\frac{13}{9}+\frac{2u_1}{9} + \frac{u_2}{9}+\frac{8}{3}b+\varepsilon}.
\end{equation}\\

\emph{Case II:} Suppose that  $|\Upsilon_1^{-1} - \alpha_1| \leq T^{-u_1}$, so that necessarily $\alpha_1 = 1$ and $\Upsilon_1^{-1} = 1 + O(T^{-u_1})$, then the coefficients simplify
\begin{displaymath}
\begin{split}
  a_9 &= (\Upsilon_1^{-1} - 1)^2\left(-\alpha_2 + O(T^{-u_1} (1 + \Upsilon_2^{-1}))\right)  \ll T^{-2u_1}, \\
  a_7 & = -C_1(\Upsilon_1^{-1} - 1)\left(-2\alpha_2 + \Upsilon_2^{-1} + O(T^{-u_1} (1 + \Upsilon_2^{-1}))\right) \ll T^{2-u_1+b},\\
  a_6 &= -C_2\left(-\Upsilon_2^{-1} + O(T^{-u_1}(1 + T^{-u_1}\Upsilon_2^{-1}))\right) \asymp  -T^3 \Upsilon_2^{-1} \gg T^{3-3b}
    \end{split}
\end{displaymath}
(using that $u_1 > b$ by \eqref{exponents}),  so
\begin{equation*}
\begin{split}
\sum_{i=0}^6 a_it^i \ll T^{\varepsilon} \left(T^{\frac{17}{2}+\frac{21}{2}b}+ |U_1|^9 T^{-2u_1}  + |U_1|^7 T^{2-u_1+b}\right) 
	\ll T^{9-u_1+8b+\varepsilon},
\end{split}
\end{equation*}
We apply the same reasoning as before to obtain
\begin{equation}\label{lem2-2}
\begin{split}
\text{{\rm meas}}(\mathcal{A}) \ll  T^{1-\frac{u_1}{6}+\frac{11}{6}b+\varepsilon} B_3 E_1^{-1/2} \leq T^{\frac{3}{2}-\frac{u_1}{6}+\frac{10}{3}b+\varepsilon} 
\end{split}
\end{equation}
in this case. \\

\emph{Case III:} Our final case is  $|\Upsilon_1^{-1}\Upsilon_2^{-1} - \alpha_2\Upsilon_1^{-1} - \alpha_1 \Upsilon_2^{-1}| \leq T^{-u_2}$, but $|\Upsilon_1^{-1} - \alpha_1| \geq T^{-u_1}$, so that  
\begin{equation}\label{ups1}
\Upsilon_1^{-1} \Upsilon_2^{-1} = \alpha_2\Upsilon_1^{-1} + \alpha_1\Upsilon_2^{-1} + O(T^{-u_2}).
\end{equation}
This implies
$\Upsilon_2 ^{-1} (\Upsilon_1^{-1} - \alpha_1) = \alpha_2 \Upsilon_1^{-1} + O(T^{-u_2})$, i.e.\ 
\begin{equation}\label{ups2}
\Upsilon_2^{-1} = \alpha_2 \Upsilon_1^{-1}(\Upsilon_1^{-1} - \alpha_1) ^{-1} + O(T^{-u_2 + u_1}).
\end{equation}
By \eqref{ups1} for the computation of $a_7$ and \eqref{ups2} for the computation of $a_6$, we then have 
 \begin{displaymath}
\begin{split}
  a_9 &\ll (\Upsilon_1^{-1} - \alpha_1)^2 T^{-u_2} \ll T^{-u_2+2b},\\
  a_7 & = -C_1\left(\alpha_1\alpha_2 \Upsilon_1^{-2} + O((1 + \Upsilon_1^{-1} + \Upsilon_1^{-2})T^{-u_2})\right) \asymp - \alpha_1\alpha_2 C_1 \Upsilon_1^{-2},\\
  a_6 & = -C_2\left( \alpha_1\alpha_2  \Upsilon_1^{-3}  (  \Upsilon^{-1}_1 - \alpha_1)^{-1} (1 - 2\alpha_1 \Upsilon_1) + O(T^{3b + u_1 - u_2})\right). 
\end{split}
\end{displaymath} 
The formula for $a_7$ requires also  $u_2 > 7b$, which is implied by \eqref{exponents}.  We now distinguish two subcases. 

\emph{Case IIIa:} Suppose first that $C_1 \geq T^{2-c}$. Then
\begin{displaymath}
\begin{split}
\sum_{i=0}^7 a_it^i \ll  T^{\varepsilon} \left( T^{\frac{17}{2}+\frac{21}{2}b}+ |U_1|^9   T^{-u_2+2b}\right) \ll 
T^{9-u_2+11b+\varepsilon}, 
\end{split}
\end{displaymath}
and again we apply the same reasoning to obtain
\begin{equation}\label{lem2-3}
\begin{split}
  \text{{\rm meas}}(\mathcal{A}) \ll T^{1-\frac{u_2}{7}+\frac{17}{7}b+\frac{c}{7} \varepsilon} B_3 E_1^{-1/2} \leq T^{\frac{3}{2}-\frac{u_2}{7}+\frac{55}{14}b + \frac{c}{7}+\varepsilon}.
 \end{split}
\end{equation}

\emph{Case IIIb:} Suppose now $C_1 \leq T^{2-c}$. We will see in a moment that  $|1 - 2\alpha_1 \Upsilon_1|  \geq T^{-d}$. Taking this for granted,  the second condition in \eqref{exponents} implies that   $a_6 \gg T^{3-9b-d}$ and
\begin{equation*}
\begin{split}
\sum_{i=0}^6 a_it^i \ll T^{\varepsilon} \left(T^{\frac{17}{2}+\frac{21}{2}b}+ |U_1|^9 T^{-u_2+2b}  + |U_1|^7 T^{2+2b -c}\right) 
	\ll T^{9-u_2+11b+\varepsilon} + T^{9 + 9b - c}.
\end{split}
\end{equation*}
Therefore, by the same reasoning as before, 
\begin{equation}\label{case3b}
\begin{split}
&\text{meas}(\mathcal{A})\\
& \ll \left(T^{1 - \frac{u_2}{6} + \frac{10}{3} b + \frac{d}{6}} + T^{1 + 3b - \frac{c}{6} + \frac{d}{6}}\right)B_3E_1^{-1/2} T^{\varepsilon} \ll T^{\frac{3}{2} - \frac{u_2}{6} + \frac{29}{6} b + \frac{d}{6}+\varepsilon} + T^{\frac{3}{2} + \frac{9}{2}b - \frac{c}{6} + \frac{d}{6}+\varepsilon}.
\end{split}
\end{equation}
It remains to show that under the current conditions the case   $|1 - 2\alpha_1 \Upsilon_1|  \leq T^{-d}$ cannot happen. Indeed, if this was the case, then
  $\alpha_1 = 1$ and $\Upsilon_1^{-1} =2 + O(T^{-d})$. Then from \eqref{ups2} we obtain  $\alpha_2 = 1$ and $\Upsilon_2^{-1}  =2+ O(T^{-d} + T^{-u_2 + u_1}).$ Substituting back into \eqref{cubic}, we obtain
\begin{displaymath}
\begin{split}
 & t_1^3 - t_1^2t_2 = 2C_2 + O\left(T^{\varepsilon} B_3|U_1|^2 E_1^{-1/2} + |U_1^3| T^{-d} + |U_1| T^{2-c}    \right),\\
  &t_2^3 - t_2^2t_1 = -2C_2 + O\left(T^{\varepsilon} B_3|U_2|^2 E_2^{-1/2} + |U_2^3| (T^{-d} + T^{-u_2 + u_1}) + |U_2| T^{2-c}    \right).
\end{split}
\end{displaymath}
Both error terms are  $O( T^{3 + \max(3b-d, 3b-u_2 + u_1, b-c)} ) $. Our conditions \eqref{exponents} imply that the exponent is strictly less than 3. Thus subtracting the two equations gives
$$(t_1 - t_2) (t_1^2 + t_2^2) \asymp 4C_2 \asymp T^3,$$
so $t_1 - t_2  \asymp T^3/(|U_1|^2 + |U_2|^2) \gg T^{1-2b}$. On the other hand, using this bound and 
 adding the two equations gives
$$T^{3-5b} \ll  (t_1 - t_2)^2 (t_1 + t_2) \ll T^{3 + \max(3b-d, 3b-u_2 + u_1, b-c)}.$$
This contradicts \eqref{exponents} for $T$ sufficiently large. \\

We now choose $$d = 9b, \quad u_1 = \frac{3}{34}-\frac{5}{2}b, \quad u_2 = \frac{13}{68}+\frac{59}{4}b, \quad c = \frac{3}{34}+\frac{27}{2}b, \quad b < \frac{3}{119}.$$
This equalizes \eqref{lem2-1}, \eqref{lem2-2}, \eqref{lem2-3} and the second term in \eqref{case3b} and satisfies \eqref{exponents}, where we need $b < 3/119$ for the first condition and $b < 7/255$ for the second. 
Combining \eqref{lem2-1}, \eqref{lem2-2}, \eqref{lem2-3}, \eqref{case3b}, we obtain
$$\text{meas}(\mathcal{A}) \ll T^{\frac{3}{2} - \frac{1}{68} +\frac{15}{4}b+\varepsilon}.$$
After substituting into \eqref{basicbound} and observing that $(B_3/(B_1B_2)) \ll T^{-1+b}$ by our assumption \eqref{case1a} and the triangle inequality \eqref{BB},  we obtain
\begin{equation}\label{lem2-all}
\mathcal{I}(B_1, B_2, B_3) \ll  \frac{\text{meas}(\mathcal{A}) }{T^{2-\varepsilon}} \left(\frac{B_3}{B_1B_2}\right)^{1/2}\ll T^{-1-\frac{1}{68} + \frac{17}{4}b + \varepsilon},
\end{equation}
provided that $b < 3/119$.

\subsection{Another special case} Here we deal with the special case 
\begin{equation}\label{case7}
|U_1| \geq T^{1+b},  \quad  T^{1- b/4} \leq  |U_2| \leq T^{1 + b/4}.
\end{equation}
We must have $B_1 \asymp B_3 \asymp |U_1|$, and by   \eqref{2b} and \eqref{basicbound} we conclude $$\mathcal{I}(B_1, B_2, B_3) \ll T^{-1+\varepsilon}.$$
Moreover, by \eqref{consist} we have 
\begin{equation}\label{special1}
\Upsilon_2 \asymp \frac{(T+ |U_2|)^2 B_2 }{|U_1| |U_2|^2} \ll \frac{(T+ |U_2|)^3   }{|U_1| |U_2|^2} \ll \frac{T^{1+b/2}}{|U_1|} \leq T^{-b/2}.
\end{equation}
Since   $E_1  \asymp |U_1|$, it follows from
$$T^{\varepsilon} |U_1|^{-1/2}   \gg g_1(t_1, t_2; d, r) =   - \log \Upsilon_1 + O\left(T^{1 + b/4 }|U_1|^{-1}\right)$$
that 
\begin{equation}\label{special2}
  |\Upsilon_1 - 1| \ll T^{1+\frac{b}{4}  } |U_1|^{-1} + T^{\varepsilon} |U_1|^{-1/2}\ll T^{-b/2}.
  \end{equation}
This suffices for the proof of \eqref{11cbound1}. Of course the same argument works with exchanged indices. 

\subsection{The remaining cases} Having the previous cases out of the way, we will show the bound
\begin{equation}\label{toshow}
  \mathcal{I}(B_1, B_2, B_3) \ll T^{-1-\frac{b}{8} +\varepsilon}
\end{equation}
in all other cases; choosing $b=\Oneb$ here and in \eqref{lem2-all} then gives \eqref{11cbound2} and completes the proof of Lemma \ref{lem6}(d).    To this end we distinguish the following principal cases
\begin{displaymath}
\begin{split}
& (1) \quad |U_1| \asymp |U_2| \asymp T, \quad\quad \,\,\,  (2) \quad  T^{1+b} \geq |U_1|  \ggg  T \asymp |U_2|, \quad \quad   (3)  \quad  |U_1|  \lll  T \asymp |U_2|,  \\
& (4) \quad  |U_1| \geq |U_2| \ggg T,  \quad\quad (5) \quad   |U_1| \ggg T \ggg  |U_2|,   \quad\quad\quad\quad\quad  (6) \quad |U_1| \leq |U_2| \lll T \end{split}
\end{displaymath}
with the understanding that those situations covered by \eqref{case1a} and \eqref{case7} and its version with exchanged indices are excluded.  By symmetry, this covers all possibilities. 

\emph{Case  1:} By \eqref{r} we conclude $B_1 \asymp B_2$ in the present case, and hence, by \eqref{BB} and in order to stay away from \eqref{case1a}, we must have $B_3 \leq T^{1-b}$. By \eqref{1b} and \eqref{basicbound}  we obtain  $ \mathcal{I}(B_1, B_2, B_3) \ll T^{-1-b +\varepsilon}$. 

\emph{Case  2:} The present assumption  together with \eqref{r} implies $B_1 \asymp B_3 \asymp |U_1|$, $B_2 \asymp T$. This case is already covered in \eqref{case1a} and currently excluded.

\emph{Case  3:} The present assumption together with \eqref{r} implies $B_1 \asymp B_2 \asymp B_3 \asymp T$. Hence in order to stay away from \eqref{case1a} we must have $|U_1| \leq T^{1-b}$. Now \eqref{1b} and \eqref{basicbound} imply $ \mathcal{I}(B_1, B_2, B_3) \ll T^{-1-2b +\varepsilon}$.

\emph{Cases  4-6}: These are handled verbatim as in \cite[Sections 15.11 - 15.13]{BB} based only on \eqref{1b} and \eqref{2b}. In all cases we confirm \eqref{toshow}.  

\section{Completion of the proof of Theorem \ref{mainthm}}

Lemmas \ref{lem1} -- \ref{lem6} are the exact analogues of \cite[Lemma 6 -- 11]{BB}, and starting from the basic inequality \eqref{basic}, the proof of Theorem \ref{mainthm} follows now verbatim as in \cite[Sections 8 -- 11]{BB}.  For convenience we indicate the key steps. We return to \eqref{basic} and apply the Kuznetsov formula \eqref{kuz} to the right hand side. We estimate each of the four resulting terms $\Delta$, $\Sigma_{4}$, $\Sigma_{5}$ and $\Sigma_6$. 

\subsection{The $\Delta$ term} By \eqref{intF} we have
$$\Delta \ll \frac{T^{3+\varepsilon}}{M^2L^2}    \sum_{j=1}^3 \sum_{\substack{\ell_1, \ell_2 \asymp L\\ \ell_1, \ell_2 \text{ prime}}}  \sum_{\substack{r_0r_1r_2 = \ell_1^j \\ s_0s_1s_2 = \ell_2^j }} \sum_{\substack{r_2n_1 \asymp M\\ s_0m_1 \asymp M\\s_1n_2 \asymp M\\ r_0m_2 \asymp M}} \delta_{\substack{m_1s_1 = m_2r_2\\ n_2s_2 = n_1 r_1}}.$$
We distinguish the cases $\ell_1 = \ell_2$ and $\ell_1 \not= \ell_2$ and estimate both contributions trivially getting
\begin{equation}\label{Delta}
  \Delta  \ll T^{3+\varepsilon}/L.
\end{equation}

\subsection{The $\Sigma_4$ term}
We have
$$\Sigma_4= \frac{1}{L^2}\sum_{\epsilon = \pm 1} \sum_{j=1}^3\sum_{\substack{\ell_1, \ell_2 \asymp L\\ \ell_1, \ell_2 \text{ prime}}} \sum_{\substack{r_0r_1r_2 = \ell_1^j \\ s_0s_1s_2 = \ell_2^j }} |\Sigma_4(r, s)|,$$
 where 
\begin{displaymath}
\begin{split}
\Sigma_4(r, s) := \frac{T^{\varepsilon}}{M^2 }   &  \sum_{m_1, m_2, n_1, n_2}  W\Bigl(\frac{r_2n_1}{M}\Bigr) W\Bigl(\frac{s_0m_1}{M}\Bigr) \overline{W\Bigl(\frac{s_1n_2}{M}\Bigr) W\Bigl(\frac{r_0m_2}{M}\Bigr)} \\
& \times 
\sum_{\substack{D, \delta\\ n_2s_2 \delta = m_2r_2 D}} \frac{\tilde{S}(-\epsilon n_1r_1, n_2s_2, m_1s_1; D, D\delta)}{D^2\delta} \Phi_{w_4}\left( \frac{\epsilon n_1n_2m_1 r_1s_1s_2}{D^2\delta} ; d \right).
\end{split}
\end{displaymath}
We use Lemma \ref{lem3} to truncate the $D, \delta$-sum at
$$ D^{2}\delta \leq \frac{M^3}{T^{3-\varepsilon}} \cdot \frac{r_1s_2}{r_2s_0}.$$
Next we apply Poisson summation in the $m_1$, $n_1$ variables. The Fourier transform of the $w_4$-Kloosterman sum $\tilde{S}$ is estimated in Lemma \ref{lem1}. As usual, the dual variables, say $x$ and $y$, can be truncated using integration by parts, and to this end we need \eqref{w42}. This effectively restricts to
\begin{displaymath}
\begin{split}
&|x| \leq X := T^{\varepsilon} \left(\frac{M^3 r_1s_2}{r_2s_0D^2 \delta}\right)^{1/3} \frac{r_2D}{M} = T^{\varepsilon}\left(\frac{r_1r_2^2s_2D}{ s_0\delta}\right)^{1/3},\\
& |y| \leq Y := T^{\varepsilon} \left(\frac{M^3 r_1s_2}{r_2s_0D^2 \delta}\right)^{1/3} \frac{s_0\delta}{M} = T^{\varepsilon}\left(\frac{r_1s_2s_0^2\delta^2}{r_2D^2}\right)^{1/3}.
\end{split}
\end{displaymath} 
We are left with bounding 
 \begin{equation*}
\begin{split}
\Sigma_4(r, s; \mu) :=  \frac{T^{3+\varepsilon}}{r_2s_0} \sum_{D^2 \delta \leq  T^{\varepsilon} \frac{M^3r_1s_2}{T^3r_2s_0}}  \frac{Dr_1s_1  }{D^2 \delta}  \sum_{\substack{ r_0m_2 \asymp M\\ s_1n_2 \asymp M\\ n_2s_2 \delta = m_2r_2 D}} \sum_{\substack{|x| \leq X, |y| \leq Y\\ (D, x) = (r_1, x)\\ (\delta, y) = (s_1, y)}} \Bigg|\tilde{\mathcal{K}}_{d, r}\left(\frac{\epsilon   n_2  r_1s_1s_2  M^2}{D^2 \delta r_2s_0 }, \frac{ xM}{r_2D}, \frac{  yM}{s_0 \delta} \right)\Bigg|
\end{split}
\end{equation*}
with $\tilde{\mathcal{K}}_{\mu}(\Xi, U, V)$ as in \eqref{deftildeK}, and the summation conditions imply that the first argument of $\tilde{\mathcal{K}}_{d, r}$ is $\gg T^{3-\varepsilon}$, while the other two arguments are 0  or 
at least $T^{1-\varepsilon}/L^j$ (which in needed for the application of Lemma \ref{lem5} in a moment). 

We start with the central Poisson term $x=y=0$, which by Lemma \ref{lem5}(b) can be bounded trivially by 
$T^{9/2 + \varepsilon} r_1 M^{-2}.$ 
Next we use Lemma \ref{lem5}(a) to show that $xy = 0$ implies $x=y=0$, up a negligible error, so that for the remaining contribution we can assume $xy \not= 0$. In this case we apply Lemma \ref{lem5}(c) and obtain altogether
 \begin{equation*}
   \Sigma_4(r, s) \ll \frac{T^{9/2 + \varepsilon} r_1}{M} +  T^{\frac{5}{2}+\varepsilon}  \left(\frac{r_1^4 s_1^2 s_2^2}{r_0^2r_2^2 s_0}\right)^{1/3}
 \end{equation*}
and in total
\begin{equation}\label{w4final}
  \Sigma_4 \ll T^{\frac{9}{2}   +\varepsilon}L^3 M^{-2} + T^{\frac{5}{2}  + \varepsilon} L^6. \end{equation}

Of course, the $\Sigma_5$ term satisfies the same bound by symmetry. 

\subsection{The $\Sigma_6$ term}
We have
\begin{displaymath}
\begin{split}
\Sigma_6 &= \frac{T^{\varepsilon}}{M^2L^2}   \sum_{j=1}^3\sum_{\substack{\ell_1, \ell_2 \asymp L\\ \ell_1, \ell_2 \text{ prime}}} \sum_{\substack{r_0r_1r_2 = \ell_1^j \\ s_0s_1s_2 = \ell_2^j }}\Biggl| \sum_{\epsilon \in \{\pm 1\}^2} \sum_{m_1, m_2, n_1, n_2}  W\Bigl(\frac{r_2n_1}{M}\Bigr) W\Bigl(\frac{s_0m_1}{M}\Bigr) \overline{W\Bigl(\frac{s_1n_2}{M}\Bigr) W\Bigl(\frac{r_0m_2}{M}\Bigr)}  \\
&  \times  \sum_{D_1, D_2} 
\frac{S(\epsilon_2 n_1r_1, \epsilon_1 m_2r_2, m_1s_1, n_2s_2; D_1, D_2)}{D_1D_2} \Phi_{w_6} \left(\left(\frac{-\epsilon_2 n_1m_1s_1r_1D_2}{D_1^2}, \frac{ -\epsilon_1n_2m_2D_1s_2r_2}{D_2^2} \right); d\right)\Biggr|.
\end{split}
\end{displaymath}
By \eqref{J1} we can truncate the $D_1, D_2$-sum at 
\begin{equation*}
D_1 \leq T^{\varepsilon} \frac{M^2}{T^2} \left(\frac{r_1^2s_1s_2}{r_0r_2 s_0^2}\right)^{1/3} \ll \Deltavar := \frac{M^2L^j}{T^{2-\varepsilon}}, \quad D_2 \leq T^{\varepsilon} \frac{M^2}{T^2} \left(\frac{r_1r_2s_2^2}{r_0^2 s_0s_1}\right)^{1/3}  \ll \Deltavar.
\end{equation*}
We apply Poisson summation to all four variables $n_1, n_2, m_1, m_2$. The dual variables, say $x_1, x_2, y_1, y_2$ can be truncated by partial  integration and \eqref{J2} at 
$$|x_1| , |x_2|, |y_1|, |y_2| \ll T^{\varepsilon} L^j(D_1 + D_2)^{1/2} =:   X.$$
Using the notation \eqref{hatS} and \eqref{defK}, we are left with bounding
\begin{displaymath}
\begin{split}
 &  \frac{T^{ \varepsilon}M^2}{L^2}   \sum_{j=1}^3 \sum_{\substack{\ell_1, \ell_2 \asymp L\\ \ell_1, \ell_2 \text{ prime}}} \sum_{\substack{r_0r_1r_2 = \ell_1^j \\ s_0s_1s_2 = \ell_2^j }}  \Biggl|     \sum_{\epsilon \in \{\pm 1\}^2} \sum_{ D_1, D_2 \leq \Deltavar}   \sum_{\substack{ |x_1|, |x_2| \leq \Xvar\\   |y_1|, |y_2| \leq  \Xvar}} \frac{\widehat{S}_{\epsilon_2r_1, s_1, \epsilon_1 r_2, s_2}(x_1, y_1, x_2, y_2; D_1, D_2)}{r_0r_2s_0s_1D_1D_2}\\
&\times  \int_{(0)} F(r)   \mathcal{K}_{d, r}\left(\frac{-\epsilon_2 M^2  s_1r_1D_2}{D_1^2 r_2s_0}, \frac{-\epsilon_1 M^2 D_1s_2r_2}{D_2^2 s_1r_0};  \frac{x_1 M}{D_1s_0}, \frac{y_1  M}{D_1r_2} ; \frac{x_2 M}{D_2r_0}, \frac{y_2 M}{D_2s_1} \right) \textbf{spec}^d(r) \, \frac{dr}{2\pi i} \Bigg|
\end{split}
\end{displaymath}
for $d \asymp T$. The summation conditions imply
$$\min\left(\frac{M}{D_1s_0}, \frac{M}{D_1 r_2}\right) \gg \frac{T^{1/2 - \varepsilon}}{L^{2j}} $$
so that the last four entries in  $\mathcal{K}_{d, r}$ are in particular $\gg T^{\varepsilon}$ if non-zero.

We use Lemma \ref{lem6}(a) to show that we need to distinguish three cases, up to symmetry and  negligible errors: the \emph{central} term $\Sigma_6^0$ where $x_1 = x_2 = y_1 = y_2 = 0$, the \emph{mixed} terms $\Sigma^{\text{mix}}_6$, where $x_1 = y_1 = 0 \not= x_2y_2$, and the \emph{generic} terms $\Sigma_6^{\text{gen}}$ where $x_1x_2y_1y_2 \not= 0$. \\

For the central term, we distinguish the cases $\ell_1 = \ell_2$ and $\ell_1 \not= \ell_2$. In the latter we use Lemma \ref{lem2}(c) and \eqref{c2} to bound this contribution by 
$\Sigma_6^0 \ll T^{3+\varepsilon} L^{-1}.$ 
Recall that \eqref{c2} was the consequence of a hidden zero of a the sum of the $(+, -)$ and the $(-, +)$ kernel  functions in the Kuznetsov formula.  
 For the contribution of $\ell_1 = \ell_2$, Lemma \ref{lem2}(d) and the simple bound \eqref{c1} suffice to obtain the same bound, so that
 \begin{equation}\label{central}
   \Sigma^{0}_6 \ll \frac{T^{3+\varepsilon}}{L}.  
\end{equation}
We observe that here we do not save in $T$, but rather in $L$, and the central Poisson term really furnishes an off-diagonal main term.\\

For the mixed terms, we apply the bounds from Lemma \ref{lem2}(a) and (b) and from Lemma \ref{lem6}(b), and obtain after direct  computation 
 \begin{equation}\label{mix}
 \begin{split}
 \Sigma^{\text{mix}}_6  & \ll   \frac{T^{3+\varepsilon} L^{31} }{M^{1/2}}  + T^{5/2 + \varepsilon} L^{23} .
\end{split}
\end{equation}\\

The generic terms lead to the most complicated analysis. From Lemma \ref{lem2}(a)  we obtain the congruences $x_1y_1 \equiv  r_1s_1 D_2 \, (\text{mod } D_1)$, $x_2y_2 \equiv  r_2s_2 D_1 \, (\text{mod } D_2)$, which we re-write as
\begin{equation}\label{cong}
 x_1y_1=   r_1s_1 D_2 + c_1D_1, \quad x_2y_2=   r_2s_2 D_1 + c_2D_2
 \end{equation}
with $c_1, c_2 \in \Bbb{Z}$. We put all variables in dyadic ranges. 
We now distinguish the three cases $c_1 = c_2 = 0$, $c_1 c_2 = 0$ but $(c_1, c_2) \not= 0$, and $c_1 c_2 \not= 0$. We call the corresponding contributions $\Sigma_6^{\text{gen}, 0}$, $\Sigma_6^{\text{gen}, \text{mix}}$ and $\Sigma_6^{\text{gen}, \ast}$, respectively. This uses the full force of difficult part Lemma \ref{lem6}(d), and the idea is that we either save a small $T$-power from \eqref{11cbound2}, or we have the extra condition \eqref{specialcase}, which shortens certain variables and gives again a saving. In very unbalanced situations this does not suffice, but then \eqref{condition-e} will be applicable, so that \eqref{cond11c} gives a saving. 

Without any difficulty we obtain
\begin{equation}\label{gen0}
 \Sigma_6^{\text{gen}, 0}  \ll T^{3 - \Eighthbfrac + \varepsilon} L^3. 
\end{equation}
Distinguishing cases as to whether \eqref{specialcase} holds or not, we obtain
\begin{equation}\label{genast}
\Sigma_6^{\text{gen}, \ast} \ll T^{3-\Eighthbfrac+\varepsilon} L^{12}. 
\end{equation}
Finally, in the situation of $
\Sigma_6^{\text{gen}, \text{mix}}$, we can assume by symmetry $c_2 = 0 \not= c_1$. Now \eqref{cong} becomes
$ x_1y_1=   r_1s_1 D_2 + c_1D_1$,  $x_2y_2=   r_2s_2 D_1$. Picking $x_1, y_1, D_2$ determines, up to a divisor function, $c_1$ and $D_1$ for fixed $r, s$. If $D_1$ is in a small dyadic range, this does not buy us anything, and it is in this situation, where  Lemma \ref{lem6}(e) is needed to obtain additional savings. Eventually we obtain
\begin{equation}\label{genmixfin}
   \Sigma^{\text{gen}, \text{mix}}_6  \ll T^{3 - \Eighthbfrac + \varepsilon} L^{12}. 
\end{equation} 
Combining \eqref{gen0}, \eqref{genast} and \eqref{genmixfin} with \eqref{central}  and \eqref{mix} we obtain finally
\begin{equation}\label{6finalbound}
  \Sigma_6 \ll T^{\varepsilon}\left(T^{3 }L^{-1}+  T^{3 - \Eighthbfrac  }L^{12}  + T^{3} L^{31} M^{-1/2} + T^{\frac{5}{2}  }L^{23}\right). 
\end{equation}

\subsection{The endgame}
Collecting the bounds \eqref{Delta}, \eqref{w4final} and \eqref{6finalbound}, we see that under the assumptions \eqref{sizeM} and \eqref{lambda} we can bound $\mathcal{L}_M(\pi_0)^2$ in \eqref{basic} by
$$\mathcal{L}_M(\pi_0)^2 \ll T^{\varepsilon} \left(T^{3 - \Eighthbfrac  + 12\lambda}  + T^{\frac{9}{4} + \frac{\eta}{2} + 31\lambda} + T^{\frac{5}{2} + 23\lambda} + T^{\frac{3}{2} + 3\lambda + 2\eta} + T^{3 - \lambda}\right),$$
and we recall the  trivial bound 
$$\mathcal{L}_M(\pi_0)^2 \ll T^{3 - 2\eta+\varepsilon},$$
see \eqref{etasmall} and \eqref{sizeM}. Now we choose $\eta = 1/100$ and $\lambda = \deltaFinal$   to complete the proof of Theorem \ref{mainthm}.

\end{document}